\documentclass[11pt,a4paper]{article}
\usepackage[T1]{fontenc}
\usepackage[utf8]{inputenc}
\usepackage{lmodern}
\usepackage[margin=28mm]{geometry}
\usepackage{amsmath,amssymb,amsthm,mathtools}
\usepackage{enumitem,microtype,xcolor}
\usepackage[colorlinks=true,linkcolor=blue!55!black,citecolor=blue!55!black,urlcolor=blue!55!black]{hyperref}
\usepackage[nameinlink,noabbrev]{cleveref}
\crefname{theorem}{theorem}{theorems}
\Crefname{theorem}{Theorem}{Theorems}
\crefname{proposition}{proposition}{propositions}
\Crefname{proposition}{Proposition}{Propositions}
\crefname{lemma}{lemma}{lemmas}
\Crefname{lemma}{Lemma}{Lemmas}
\crefname{corollary}{corollary}{corollaries}
\Crefname{corollary}{Corollary}{Corollaries}
\crefname{definition}{definition}{definitions}
\Crefname{definition}{Definition}{Definitions}
\crefname{assumption}{assumption}{assumptions}
\Crefname{assumption}{Assumption}{Assumptions}
\crefname{remark}{remark}{remarks}
\Crefname{remark}{Remark}{Remarks}
\newtheorem{theorem}{Theorem}[section]
\newtheorem{proposition}[theorem]{Proposition}
\newtheorem{lemma}[theorem]{Lemma}
\newtheorem{corollary}[theorem]{Corollary}
\theoremstyle{definition}
\newtheorem{assumption}[theorem]{Assumption}
\newtheorem{definition}[theorem]{Definition}
\theoremstyle{remark}
\newtheorem{remark}[theorem]{Remark}
\numberwithin{equation}{section}
\allowdisplaybreaks[2]
\newcommand{\R}{\mathbb R}
\newcommand{\dd}{\,\mathrm d}
\newcommand{\Sch}{\mathcal S}
\newcommand{\divv}{\operatorname{div}}
\newcommand{\supp}{\operatorname{supp}}
\newcommand{\ad}{\operatorname{ad}}

\newcommand{\Pp}{\mathcal P}
\title{Regularization and Full Support\\
for the CBO Equation with Measure Initial Data}
\author{LUKANG SUN\\
\small Department of Mathematics, National University of Singapore\\
\small \href{mailto:lukang.sun@nus.edu.sg}{\texttt{lukang.sun@nus.edu.sg}}}
\date{October 22, 2026}

\begin{document}
\maketitle
\begin{abstract}
We prove global well-posedness and regularization for the isotropic
consensus-based optimization equation with initial probability measures
of finite fourth moment. Under suitable weight decay assumptions, after any nonconstant consensus motion,
both the density and the consensus trajectory become smooth, and all
spatial Sobolev norms are locally bounded in time. The proof uses a
Fourier commutator estimate, localization and iteration, followed by a
nonlinear time bootstrap. Counterexamples show that smoothing can fail
at a stationary consensus point. For non-Dirac initial data, we also prove
positivity away from the consensus point and full support at every
positive time in all dimensions.
\end{abstract}



\tableofcontents

\section{Introduction, main results and assumptions}

Fix $d\geq1$, $\lambda,\sigma,\alpha>0$, and $T>0$.
We study the mean-field equation of consensus-based optimization (CBO),
\begin{equation}\label{eq:cbo}
\partial_t\rho
 =\lambda\divv_v((v-m(t))\rho)
       +\frac{\sigma^2}{2}\Delta_v(|v-m(t)|^2\rho),
\qquad
m(t)=\frac{A(t)}{B(t)},
\quad w=e^{-\alpha f}.
\end{equation}
Here $\rho_t$ is a probability measure, $f:\R^d\to\R$ is the objective
function, $A(t)=\int v w(v)\dd\rho_t(v)$, and
$B(t)=\int w(v)\dd\rho_t(v)$. Thus $m(t)$ is the barycenter of
$\rho_t$ weighted by $e^{-\alpha f}$. The structure of \eqref{eq:cbo},
an affine drift together with a diffusion that vanishes quadratically
at a single point, is shared by other equations, among them the
Fokker--Planck equations of wealth distribution
\cite{TorregrossaToscani2018Debts}, of quickest change detection
\cite{Peskir2006}, and of neuronal membrane potentials
\cite{DOnofrioLanskyPirozzi2018}; our results apply to these as well,
see \cref{rem:affine-applications}.

Consensus-based optimization was introduced by Pinnau, Totzeck, Tse,
and Martin \cite{PinnauTotzeckTseMartin2017} as a derivative-free
method for global optimization: particles are drawn toward a weighted
average of their positions and explore by a noise that switches off
as they reach it. The idea has proved remarkably flexible. Within a few
years it was taken up in numerical analysis
\cite{HaJinKim2021,GrassiPareschi2021,CiprianiHuangQiu2022,GerberHoffmannVaes2025},
nonconvex optimization \cite{CCTT,FKR2024,FornasierSunMultiple2025},
machine learning
\cite{CarrilloJinLiZhu2021,FornasierHuangPareschiSunnen2021,CarrilloGarciaTrillosLiZhu2024},
sampling and inverse problems
\cite{CarrilloHoffmannStuartVaes2022,BungertRoithWacker2025,BouillonEtAl2026},
and finance \cite{BaeEtAl2022}, and the same mechanism was carried
over to particles constrained to hypersurfaces
\cite{FornasierHuangPareschiSunnen2020} or to general feasible sets
\cite{BorghiHertyPareschi2023Constraints,BeddrichEtAl2026}, to
objectives with several global minimizers
\cite{BungertRoithWacker2025,FornasierSunMultiple2025}, to Pareto
fronts \cite{BorghiHertyPareschi2023Multi,KlamrothStiglmayrTotzeck2024},
to objectives given by expectations \cite{BonandinHerty2025}, and to
min--max problems \cite{HuangQiuRiedl2024}, multiplayer games
\cite{ChencheneHuangQiu2025}, and bilevel problems
\cite{GarciaTrillosLiRiedlZhu2024}. The dynamics themselves have been
varied many ways: componentwise noise
\cite{CarrilloJinLiZhu2021,FornasierKlockRiedl2022}, a sampling variant
that targets a distribution rather than a point
\cite{CarrilloHoffmannStuartVaes2022,BouillonEtAl2026}, personal-best
memory \cite{TotzeckWolfram2020}, gradient information \cite{Riedl2024},
adaptive momentum \cite{ChenJinLyu2022}, second-order dynamics
\cite{GrassiPareschi2021,CiprianiHuangQiu2022}, mirror-descent
geometry \cite{BungertEtAl2025Mirror}, random batches
\cite{KoHaJinKim2022}, truncated noise
\cite{FornasierRichtarikRiedlSun2025}, jumps
\cite{KaliseSharmaTretyakov2023}, and feedback control
\cite{HuangHertyKaliseKantas2026}. The literature is by now too large
to survey here, and we cite only a selection.

What makes \eqref{eq:cbo} hard to analyze is visible in the equation
itself. The consensus point depends on the whole distribution, which
already complicates existence and the study of concentration
\cite{PinnauTotzeckTseMartin2017,CCTT}. More seriously for regularity,
the diffusion vanishes exactly at $v=m(t)$ and grows quadratically at
infinity, so standard parabolic theory does not apply as it stands
\cite{WangLiHuang2025}. Near the consensus point the drift is weak as
well, and nothing in the equation forces smoothness or positivity
there: whether they hold \emph{depends on how drift and diffusion act together}, not on the mere presence of noise \cite{FS}.  The
PDE theory is accordingly still incomplete where the initial datum is
rough. Probability-measure solutions exist and are unique under moment
and growth assumptions on the objective \cite{CCTT,GerberHoffmannVaes2025};
for initial densities in $L^1\cap L^\infty$, weak solutions can be built,
and $H^2$ regularity is propagated from $H^2$ initial data
\cite{WangLiHuang2025}; if the initial density is smooth, the solution
is classical and stays smooth \cite{FS}. Owing to the degenerate and nonlocal structure of \eqref{eq:cbo},
no result is available so far on the regularization of solutions
starting from a singular probability measure, that is, on whether, and
from which time on, such a solution admits a smooth density.

Our first purpose is to understand how consensus motion affects
regularization. A stationary consensus point can preserve singularities. We show that, once the
consensus point starts moving, the solution becomes smooth in space and
time and has finite spatial $H^s$ norms for every $s\geq0$
(\cref{thm:main}). This holds even for a $C^2$ objective satisfying
our weight assumptions, without an initial density or Sobolev regularity.

Our second purpose is to establish full support at positive times
from non-Dirac initial measures. This property is relevant to global
optimization: the convergence theorem of Fornasier--Klock--Riedl
\cite{FKR2024} requires initial mass in every neighborhood
of a global minimizer. Fornasier--Sun \cite{FS} established
the positivity mechanism for classical solutions starting from initial
densities in $H^m$. We extend this mechanism to measure initial data
in every dimension, including dimension one
(\cref{thm:fs-positivity}).

\begin{definition}[End of the initial stationary interval]
\label[definition]{def:onset}
For a continuous consensus trajectory $m:[0,T]\to\R^d$, set
\(t_c:=\sup\{t\in[0,T]:m(s)=m(0)\text{ for all }s\in[0,t]\}\).
\end{definition}
The defining set contains zero and is closed and downward
closed. Hence $m=m(0)$ on $[0,t_c]$. If $t_c<T$, then
$m$ is nonconstant on $[0,b]$ for every $b\in(t_c,T]$.
The cases $t_c=0$ and $t_c=T$ are allowed; the latter
means that $m$ is constant on the whole interval.

\begin{theorem}[Global well-posedness and regularity relative to consensus motion]
\label[theorem]{thm:main}
Under \cref{ass:main}, there exists a unique global weak solution in
the sense of \cref{def:weak-solution} with initial value $\rho_0$.
On each finite interval $[0,T]$, the solution is $W_2$ continuous
and has bounded fourth moments. The following regularity conclusions
hold for every finite $T>0$. Let $t_c$ be given by \cref{def:onset}
for its consensus trajectory. Then
\begin{equation}\label{eq:consensus-onset-regularity}
m\in C^k([0,T];\R^d)
       \cap C^\infty([0,T]\setminus\{t_c\};\R^d),
\end{equation}
where derivatives at endpoints are one-sided.
If $t_c=T$, the trajectory is constant, hence smooth
on the entire closed interval as well.
The given measure solution has the following regularity.
\begin{enumerate}[label=\textup{(\roman*)},leftmargin=*]
\item \emph{Initial stationary phase.} If $t_c>0$, then
\(\rho\in C^\infty\big((\R^d\setminus\{m(0)\})\times(0,t_c)\big)\).
No regularity at the stationary center is asserted.
\item \emph{After the consensus point starts moving.} If $t_c<T$, then
\begin{equation}\label{eq:global-after-motion}
\begin{aligned}
&\rho\in C^\infty(\R^d\times(t_c,T)),\\
&\rho\in C((t_c,T];H^s(\R^d))
 \quad\text{for every }s\geq0.
\end{aligned}
\end{equation}
In particular, every such Sobolev norm is bounded on
compact time intervals strictly after $t_c$.
\item \emph{An interior transition time.} If $0<t_c<T$,
then for every $v_0\neq m(t_c)$ there is a neighborhood
of $(v_0,t_c)$ on which $\rho\in C^{\infty,k+1}$.
Moreover,
\(m^{(j)}(t_c)=0, 1\leq j\leq k\).
Here $C^{\infty,k+1}$ means joint continuity of every
$\partial_v^\beta\partial_t^j$ with $|\beta|<\infty$
and $0\leq j\leq k+1$.
\end{enumerate}
\end{theorem}

\paragraph{Weak probability-measure solution}
Carrillo--Choi--Totzeck--Tse
\cite{CCTT}
prove well-posedness by a Leray--Schauder argument for bounded objectives
or objectives with two-sided quadratic growth, subject to their
increment condition. Gerber--Hoffmann--Vaes
\cite{GerberHoffmannVaes2025}
allow general two-sided polynomial growth, together with a polynomially
weighted local Lipschitz condition. 
Our weight conditions $w>0$ and $w,v_iw\in\mathcal X^2$ also cover
unbounded objectives growing more slowly than every positive power.
For example, $f(v)=[\log(1+|v|^2)]^2$ satisfies \cref{ass:main}(i)
for every $k\geq1$ and $\alpha>0$. Every positive-order derivative of
this $f$ is bounded, so $|\partial^\gamma w(v)|\leq C_{\alpha,\gamma}w(v)$.
Since $(1+|v|^2)^Nw(v)\leq\exp(N^2/(4\alpha))$ for every $N\geq0$,
we have $w,v_iw\in\Sch$.
On the other hand, $f(v)\to\infty$ and $f(v)/|v|^p\to0$ as
$|v|\to\infty$ for every $p>0$.
Thus this example satisfies neither the bounded-objective alternative
nor the positive-power polynomial lower bounds required in the
nonlinear well-posedness results cited above.

Our strategy is to first construct a local weak solution and then prove
that it can be continued to arbitrary finite times.
For local existence, we truncate the consensus denominator and use
uniform moment bounds and compactness to pass to the limit in a
time-discrete approximation
(\cref{lem:cutoff-existence,thm:local-existence}).
The key difficulty in global continuation is preventing finite-time
blow-up of $m$: the transition-law lower bound in \cref{lem:global-lower},
together with rapid decay of $w$, controls the consensus denominator and
rules out such blow-up, yielding \cref{thm:global-existence}.
The $W_2$ stability estimate in \cref{prop:weak-uniqueness} gives
uniqueness independently of the regularization argument.

\paragraph{H\"ormander condition.}
The starting point of this work is our observation that motion of the
consensus point can restore H\"ormander's bracket condition at the degenerate point.
With $y=v-m(t)$, expanding \eqref{eq:cbo} gives
\begin{equation}\label{eq:forward}
\partial_t\rho
 =\frac{\sigma^2}{2}|y|^2\Delta_v\rho
     +(\lambda+2\sigma^2)y\cdot\nabla_v\rho
     +d(\lambda+\sigma^2)\rho.
\end{equation}
For the classical bracket calculation, assume $m\in C^\infty$ and set
$X_{ij}=\frac{\sigma}{\sqrt2}y_j\partial_{v_i}$ for $1\leq i,j\leq d$
and $X_0=\partial_t-(\lambda+\frac{3\sigma^2}{2})y\cdot\nabla_v$.
Since $\sum_{i,j}X_{ij}^2=\frac{\sigma^2}{2}|y|^2\Delta_v
+\frac{\sigma^2}{2}y\cdot\nabla_v$, equation \eqref{eq:forward} becomes
$(\sum_{i,j}X_{ij}^2-X_0+d(\lambda+\sigma^2))\rho=0$.
All these fields are smooth at $y=0$. Using $[Y,Z]=YZ-ZY$, we have
$[y\cdot\nabla_v,y_j\partial_{v_i}]=0$ and
$[\partial_t,y_j\partial_{v_i}]=-m_j'(t)\partial_{v_i}$, hence
$[X_0,X_{ij}]=-\frac{\sigma}{\sqrt2}m_j'(t)\partial_{v_i}$.
For $\ell\geq1$, direct differentiation gives
$[X_0,m_j^{(\ell)}(t)\partial_{v_i}]
=(m_j^{(\ell+1)}(t)+(\lambda+\frac{3\sigma^2}{2})
m_j^{(\ell)}(t))\partial_{v_i}$.
Induction therefore yields
$\ad_{X_0}^{\,r}X_{ij}
=-\frac{\sigma}{\sqrt2}\sum_{\ell=0}^{r-1}\binom{r-1}{\ell}
(\lambda+\frac{3\sigma^2}{2})^{r-1-\ell}
m_j^{(\ell+1)}(t)\partial_{v_i}$ for $r\geq1$, where
$\ad_{X_0}Z=[X_0,Z]$. Away from $v=m(t)$, the fields $X_{ij}$ span all spatial directions.
At $v=m(t)$, let $r\geq1$ be the first index with $m^{(r)}(t)\neq0$
and choose $j$ with $m_j^{(r)}(t)\neq0$. Then
$\ad_{X_0}^{\,r}X_{ij}
=-\frac{\sigma}{\sqrt2}m_j^{(r)}(t)\partial_{v_i}$ for $i=1,\ldots,d$.
Together with the time direction in $X_0$, these vectors span all
$d+1$ space--time directions. H\"ormander's theorem
\cite{Hormander1967} thus gives local $C^\infty$
regularity wherever this condition holds.

\paragraph{$H^s$ estimate and Bootstrapping.}
The classical H\"ormander theorem
\cite{Hormander1967} requires $C^\infty$
coefficients, hence $m$ needs to be smooth.
Its local hypoellipticity conclusion does not directly give the global
$H^s$ estimates in \eqref{eq:global-after-motion}.
Our Fourier estimates require only $m\in C^1$, which follows from
\cref{prop:observables}(ii).
Any nonconstant motion then provides an interval on which one component
of $m'$ is bounded away from zero. On this interval,
\cref{lem:annular-subelliptic} gives an $R^{2/3}$ gain on frequency
annuli without differentiating $m'$ in time.
\Cref{lem:localized-annular-gain} localizes this estimate, and
\cref{prop:fourier-smoothing} iterates it to obtain all spatial $H^s$
orders. The weighted Fourier energy estimate in
\cref{lem:fourier-propagation} propagates this regularity to every later
time. The PDE and the consensus formula then bootstrap time regularity.
The proof of \cref{prop:earlier-motion} shows that $m\in C^r$, $r\geq1$,
gives $\rho\in C^{\infty,r+1}$ through the PDE.
In \cref{prop:post-motion-bootstrap}, the spatial Sobolev bounds and rapid
decay of $w$ justify differentiating the weighted integrals in the
consensus formula, giving $m\in C^{r+1}$ without taking further
derivatives of $w$. Iterating this step makes both $m$ and $\rho$ smooth
to all orders after $t_c$, even under the order-two weight assumption.

\begin{remark}[Sharpness at a stationary consensus point]
\label[remark]{rem:stationary-sharpness}
The exclusion of the stationary center in \cref{thm:main}(i) is essential.
The simplest example is $\rho_0=\delta_{v_0}$: the solution is
$\rho_t=\delta_{v_0}$ with $m(t)\equiv v_0$, so the singularity persists
at the consensus point for all time. This obstruction also occurs when
the initial measure has a density. For $f(v)=|v|^2$, \cref{sec:counterexample} constructs
global weak solutions with $m(t)\equiv0$ and the following properties:
\begin{enumerate}[label=\textup{(\roman*)},leftmargin=*]
\item A bounded, continuous, strictly positive initial density with all
polynomial moments gives a solution that is smooth away from the origin
but fails to be $C^1$ there at every positive time
(\cref{prop:counter-bounded}).
\item For every $s_0\geq0$ and $\varepsilon>0$, there is an initial
probability density with full support and all polynomial moments such
that $\rho_0,\rho_t\in H^{s_0}(\R^d)$ for every $t>0$, but $\rho_t$
has no $H^{s_0+\varepsilon}$ regularity on any neighborhood of the origin
(\cref{thm:counter-no-gain}).
\end{enumerate}
The quadratic objective and both density examples satisfy
\cref{ass:main} for every $k\geq1$. Thus neither smoothing at a
stationary center nor any positive unweighted Sobolev gain is guaranteed,
even for initial densities.
The stationary consensus in these density examples is enforced by the
shared symmetry of $f$ and $\rho_0$. For generic choices of $f$ and
non-atomic $\rho_0$, we expect the consensus to move immediately,
so that $t_c=0$.
\end{remark}

\begin{remark}[Other equations with quadratic diffusion degeneracy]
\label[remark]{rem:affine-applications}
Equation \eqref{eq:cbo} belongs to the family
\begin{equation}\label{eq:affine-family}
\partial_t\rho
 =\divv_v\bigl((\kappa(t)v-\beta(t))\rho\bigr)
   +\frac{\sigma^2}{2}\Delta_v\bigl(|v-a(t)|^2\rho\bigr),
\qquad \beta(t),a(t)\in\R^d,\ \kappa(t)\in\R,
\end{equation}
with $\beta=\lambda m$, $\kappa=\lambda$ and $a=m$. Several equations in
the literature are of this form with $d=1$. The wealth distribution
equation of Torregrossa--Toscani \cite{TorregrossaToscani2018Debts}
reads $\partial_t\rho=\lambda\partial_v((v-1)\rho)
+(\sigma^2/2)\partial_{vv}(v^2\rho)$; the Kolmogorov--Shiryaev
equation of quickest change detection \cite{Peskir2006} reads
$\partial_t\rho=-\partial_v((1+\mu v)\rho)
+(\sigma^2/2)\partial_{vv}(v^2\rho)$; and the membrane potential of a
neuron driven by an inhomogeneous geometric Brownian motion (IGBM)
satisfies $\partial_t\rho=\partial_v((v/\theta-\mu)\rho)
+(\sigma^2/2)\partial_{vv}((v-V_I)^2\rho)$, where $\theta>0$ is the
membrane time constant, $\mu$ the input and $V_I$ the inhibitory
reversal potential \cite{DOnofrioLanskyPirozzi2018}. In all of them
the diffusion vanishes at $v=a(t)$, and this point is known to cause
difficulties: concentration at zero wealth is hard to control
\cite{TorregrossaToscani2018Debts}, and the Kolmogorov--Shiryaev
equation requires a separate analysis at the boundary point zero
\cite{Peskir2006}.

The $H^s$ estimates of \cref{thm:main} extend to
\eqref{eq:affine-family} with prescribed continuous $\beta,\kappa$ and
$a\in C^1$ after a change of variables. Choose a constant
$\lambda>0$ as the contraction rate of the transformed CBO drift.
For the Kolmogorov--Shiryaev and IGBM equations, this rate can
be chosen arbitrarily and is absent from the original models.
For this fixed $\lambda$, define
$r(t)=\exp\bigl(\int_0^t(\kappa(s)-\lambda)\dd s\bigr)$ and
$q(t)=\int_0^tr(s)(\beta(s)-\lambda a(s))\dd s$. In the variable
$y=r(t)v-q(t)$ the drift is $-\lambda(y-\widetilde m(t))$ and the
diffusion matrix is $\sigma^2|y-\widetilde m(t)|^2I_d$, where
$\widetilde m(t)=r(t)a(t)-q(t)$; that is,
$\psi_t=(v\mapsto r(t)v-q(t))_\#\rho_t$ solves the linear equation in
\eqref{eq:cbo} with center $\widetilde m(t)$. Moreover
$
\widetilde m'(t)=r(t)\bigl(a'(t)-\beta(t)+\kappa(t)a(t)\bigr),
$
hence what matters is not whether the degeneracy point $a(t)$ moves,
but whether its velocity $a'(t)$ differs from the drift
$\beta(t)-\kappa(t)a(t)$ at that point. In the three examples this is
the case for all $t$. For the wealth distribution equation, we choose
the target contraction rate equal to the original relaxation
coefficient $\lambda$. Then $(\beta,\kappa,a)=(\lambda,\lambda,0)$,
hence $r=1$, $y=v-\lambda t$,
$\widetilde m(t)=-\lambda t$ and $\widetilde m'(t)=-\lambda\neq0$. For the Kolmogorov--Shiryaev
equation one has $(\beta,\kappa,a)=(1,-\mu,0)$, hence
$r(t)=e^{-(\mu+\lambda)t}$, $y=r(t)v-\int_0^tr(s)\dd s$ and
$\widetilde m'(t)=-r(t)$. For the IGBM equation one has
$(\beta,\kappa,a)=(\mu,\theta^{-1},V_I)$, hence
$\widetilde m'(t)=r(t)(V_I/\theta-\mu)$, which is nonzero in the
neuronal parameter regime $\mu>V_I/\theta$. In all three cases
$\widetilde m$ is smooth and $\widetilde m'$ has no zeros, so
\cref{prop:earlier-motion} applies on every time interval and gives,
for every initial probability measure $\rho_0\in\Pp_4(\R)$, that
$\psi\in C^\infty(\R\times(0,\infty))$ and
$\psi\in C([\tau,T];H^s(\R))$ for all $s\geq0$ and $0<\tau<T<\infty$.
Since $\rho_t$ is the image of $\psi_t$ under $y\mapsto(y+q(t))/r(t)$,
an affine map with $r(t)>0$ depending smoothly on $t$, the same holds
for $\rho$.
As far as we know, in the literature, no result of this type seems to be known for these equations:
Torregrossa--Toscani \cite{TorregrossaToscani2018Debts} prove that
the homogeneous Sobolev norm of the solution stays finite if it is
finite initially, Peskir \cite{Peskir2006} obtains an integral
representation of the fundamental solution of the Kolmogorov--Shiryaev
equation, and D'Onofrio, Lansky and Pirozzi
\cite{DOnofrioLanskyPirozzi2018} study the distribution and the
first-passage times of the IGBM.
\end{remark}

\paragraph{Related literature.}
The theoretical analysis connects particle dynamics, mean-field limits,
and convergence to global minimizers.
Carrillo--Choi--Totzeck--Tse \cite{CCTT}
studied consensus formation and approximate global optimality of the
mean-field dynamics.
Fornasier--Klock--Riedl
\cite{FKR2024} established quantitative
finite-accuracy global convergence in mean-field law using a
nonasymptotic Laplace principle and an initial support condition.
Fornasier--Sun \cite{FornasierSunMultiple2025} analyzed
a CBO variant for objectives with multiple global minimizers and proved
concentration toward the minimizing set.

Of particular relevance here is the regularity theory of Fornasier--Sun
\cite{FS}, which constructs classical
solutions of degenerate drift--diffusion equations by Galerkin
approximation and weighted Sobolev estimates for regular initial data.
Their off-center positivity result for $d>1$
\cite{FS} generates full support at positive
times and connects this PDE analysis to the global convergence mechanism
of \cite{FKR2024}. The present work studies instantaneous
regularization from measure initial data, identifies consensus motion as
the mechanism removing the degeneracy, and extends the positivity
argument to dimension one.

For background on hypoelliptic regularity,
we recall the following works. Bouchut \cite{Bouchut2002}
developed Fourier and commutator methods that transfer regularity through
transport. Alexandre
\cite{Alexandre2012}
gave energy proofs of fractional kinetic estimates. Lerner, Morimoto, and Pravda-Starov
\cite{LMP2012}
obtained sharp parameter-dependent estimates and local Sobolev
regularity using Wick quantization. Morimoto--Xu
\cite{MorimotoXu2009} proved ultra-analytic
smoothing for a linear kinetic Fokker--Planck equation using exponential
Fourier weights. These works provide background on hypoellipticity, but their derivations
do not apply directly to the CBO Fourier operator, whose structure differs
from the kinetic operators above. We independently derive estimates for
the scaled $F_R(t,\xi)=F(t,R\xi)$ on fixed frequency annuli,
with constants uniform in $R\geq1$. In particular,
\cref{lem:annular-subelliptic} uses the skew-adjoint operator
$\mathcal T_R$ to obtain a new commutator estimate assuming only $m\in C^1$.

\paragraph{Assumptions and solution convention.}

Throughout the rest of the paper, $C$ denotes a positive constant whose value may change from
line to line. When a lemma or proposition says that $C$ depends only on
certain quantities, this is meant literally: $C$ does not depend on the
frequency scale $R$, on the function being estimated, or on the solution,
except through the quantities listed.

For an integer $r\geq0$ write
\[
\mathcal X^r
 =\left\{\psi\in C^r(\R^d):
  \sup_v(1+|v|)^N|\partial^\gamma\psi(v)|<\infty
  \text{ for every }N\geq0,\ |\gamma|\leq r\right\}.
\]
The derivatives through order $r$ decay faster than every
polynomial.
\begin{assumption}[Objective and initial data]
\label[assumption]{ass:main}
Fix an integer $k\geq1$ and impose the following conditions.
\begin{enumerate}[label=\textup{(\roman*)},leftmargin=*]
\item \emph{Finite-order objective and decay of the weight.}
The objective satisfies $f\in C^{2k}(\R^d)$.
For $w=e^{-\alpha f}>0$ and every
$\psi\in\{w,v_1w,\ldots,v_dw\}$,
\(\sup_{v\in\R^d}(1+|v|)^N|\partial^\gamma \psi(v)|<\infty \quad\text{for every }N\geq0\text{ and }|\gamma|\leq2k\);
equivalently, $w,v_iw\in\mathcal X^{2k}$.
\item \emph{Initial probability measure.}
The initial datum is $\rho_0\in\Pp_4(\R^d)$, that is,
\(\rho_0(\R^d)=1, \rho_0\geq0, \int_{\R^d}|v|^4\dd\rho_0(v)<\infty\).
No density or initial Sobolev regularity is required.
\end{enumerate}
\end{assumption}

A concrete sufficient condition for \cref{ass:main}(i),
expressed entirely in terms of the objective, is
\begin{equation}\label{eq:objective-sufficient}
\begin{gathered}
f\in C^{2k}(\R^d),\qquad f(v)\geq c|v|^2-C,
\quad c>0,\\
|\partial^\gamma f(v)|
 \leq C_\gamma(1+|v|)^{N_\gamma},
\qquad 1\leq|\gamma|\leq2k.
\end{gathered}
\end{equation}
Here $c,C$ and $C_\gamma,N_\gamma$ are fixed objective
bounds: they may depend on $f,d$ and, for the indexed
constants, $\gamma$, but not on $v,t$ or the solution.

\begin{definition}[Weak probability-measure solution]
\label[definition]{def:weak-solution}
Let $T>0$. For a trajectory $m$, define the backward test operator by
\begin{equation}\label{eq:generator}
L_{m(t)}\phi(v)
 =-\lambda(v-m(t))\cdot\nabla_v\phi(v)
       +\frac{\sigma^2}{2}|v-m(t)|^2\Delta_v\phi(v).
\end{equation}
A family of probability measures
$(\rho_t)_{t\in[0,T]}$ on $\R^d$ is a weak solution of
\eqref{eq:cbo} with initial value $\rho_0$ if:
\begin{enumerate}[label=\textup{(\roman*)},leftmargin=*]
\item The curve is narrowly continuous on $[0,T]$:
for every bounded continuous $\psi:\R^d\to\R$,
the function $t\mapsto\int\psi\dd\rho_t$ is continuous.
The curve takes the prescribed value $\rho_0$ at $t=0$;
hence $\rho_t\rightharpoonup\rho_0$ as
$t\downarrow0$.
\item With $m(t)$ given by the consensus formula in
\eqref{eq:cbo}, the integration-by-parts identity
\begin{equation}\label{eq:weak}
\int_0^T\int_{\R^d}
\big(\partial_t\varphi+L_{m(t)}\varphi\big)\dd\rho_t(v)\dd t=0
\end{equation}
holds for every
$\varphi\in C_c^\infty(\R^d\times(0,T))$.
\end{enumerate}
\end{definition}

\paragraph{Arrangement of the paper.}
The proof of \cref{thm:main} is given in
Sections~\ref{sec:preliminary-regularity}--\ref{sec:finite-full} and
Appendices~\ref{sec:existence}--\ref{sec:moment-fourier-proofs}.
Some proofs involve lengthy verifications and are deferred to the
appendices to keep the main exposition uninterrupted.
\Cref{sec:preliminary-regularity} establishes the time regularity of $m(t)$
and local regularity of $\rho$ away from the consensus point.
\Cref{sec:finite-full} develops the Fourier estimates and nonlinear
bootstrap.
\Cref{sec:positivity} proves positivity away from the consensus point
and full support at positive times.
Appendix~\ref{sec:counterexample} constructs the stationary counterexamples,
and Appendix~\ref{sec:existence} proves global well-posedness.
Appendix~\ref{sec:moment-fourier-proofs} contains the proofs of the
moment bound and Fourier equation used in \cref{sec:finite-full}.

\section{Preliminary regularity}
\label{sec:preliminary-regularity}

\subsection{Properties of \texorpdfstring{$m(t)$}{m(t)}}
\label{sec:observables}

We first observe that $L_m:\mathcal X^r\to\mathcal X^{r-2}$
for $r\geq2$.

\begin{lemma}\label[lemma]{lem:observable}
If $m$ is continuous and $\phi\in\mathcal X^2$, then
$\Phi_\phi(t):=\int\phi\dd\rho_t$ is $C^1$ and
\begin{equation}\label{eq:observable}
\Phi_\phi'(t)=\int L_{m(t)}\phi\dd\rho_t.
\end{equation}
At endpoints the assertion holds with one-sided derivatives.
\end{lemma}

\begin{proof}
Let $\chi\in C_c^\infty(B_2)$ satisfy $0\leq\chi\leq1$ and
$\chi=1$ on $B_1$, and put $\chi_R(v)=\chi(v/R)$.
Choose $\phi_n\in C_c^\infty(B_{3R})$ with
$\phi_n\to\chi_R\phi$ in $C^2$.
Fix $0<s<t<T$. Choose a nonnegative $\zeta\in C_c^\infty((-1,1))$
with $\int\zeta=1$, and put
$\zeta_\varepsilon(\tau)=\varepsilon^{-1}\zeta(\tau/\varepsilon)$.
For $0<\varepsilon<\min\{s,T-t,(t-s)/2\}$, take
$\eta_\varepsilon(\tau)=\int_s^t\zeta_\varepsilon(\tau-r)\dd r$.
This function equals one on $[s+\varepsilon,t-\varepsilon]$,
vanishes outside $(s-\varepsilon,t+\varepsilon)$, and satisfies
$\eta_\varepsilon'(\tau)=\zeta_\varepsilon(\tau-s)-\zeta_\varepsilon(\tau-t)$.
Use $\varphi_n(v,\tau)=\phi_n(v)\eta_\varepsilon(\tau)$ in
\eqref{eq:weak} and let $n\to\infty$. This gives
\(\int_0^T\eta_\varepsilon'(\tau)\int\chi_R\phi\dd\rho_\tau\dd\tau
+\int_0^T\eta_\varepsilon(\tau)
\int L_{m(\tau)}(\chi_R\phi)\dd\rho_\tau\dd\tau=0\).
As $\varepsilon\to0$, narrow continuity makes the first term converge to
$\int\chi_R\phi\dd\rho_s-\int\chi_R\phi\dd\rho_t$.
The second term converges by dominated convergence to
$\int_s^t\int L_{m(\tau)}(\chi_R\phi)\dd\rho_\tau\dd\tau$.
Consequently,
\begin{equation}\label{eq:observable-cutoff-identity}
\int\chi_R\phi\dd\rho_t-\int\chi_R\phi\dd\rho_s
 =\int_s^t\int L_{m(\tau)}(\chi_R\phi)\dd\rho_\tau\dd\tau.
\end{equation}
Narrow continuity and boundedness of the time integrand extend this
identity to $0\leq s\leq t\leq T$.

To remove the cutoff, write $y=v-m(\tau)$. The product rule gives
\begin{equation}\label{eq:cutoff}
L_m(\chi_R\phi)-\chi_RL_m\phi
=-\lambda(y\cdot\nabla\chi_R)\phi+\frac{\sigma^2}{2}|y|^2
       \big((\Delta\chi_R)\phi+2\nabla\chi_R\cdot\nabla\phi\big).
\end{equation}
This commutator vanishes outside $R\leq|v|\leq2R$ and is bounded
there by $C(|\phi|+R|\nabla\phi|)$, uniformly in $\tau\in[0,T]$.
The difference $L_m(\chi_R\phi)-L_m\phi$ is the sum of this
commutator and $(\chi_R-1)L_m\phi$; the latter vanishes on $B_R$
and is bounded by $C(1+|v|)^2(|\nabla\phi|+|\nabla^2\phi|)$.
Since $\phi\in\mathcal X^2$, for every $N\geq0$ we obtain
\(\|\chi_R\phi-\phi\|_\infty
+\sup_{\tau\in[0,T]}\|L_{m(\tau)}(\chi_R\phi)-L_{m(\tau)}\phi\|_\infty
\leq C_{\phi,N}R^{-N}\).
Here $C$ depends only on $d,\lambda,\sigma,\sup_{[0,T]}|m|$ and
$\chi$; $C_{\phi,N}$ additionally depends on $N$ and the weighted
supremum norms of $\phi$ through order two with weight
$(1+|v|)^{N+2}$.

Because every $\rho_\tau$ has mass one,
\(\left|\int(\chi_R-1)\phi\dd\rho_t\right|
\leq C_{\phi,N}R^{-N}\), with the same bound at $s$, and
\(\left|\int_s^t\int
[L_{m(\tau)}(\chi_R\phi)-L_{m(\tau)}\phi]\dd\rho_\tau\dd\tau\right|
\leq C_{\phi,N}(t-s)R^{-N}\).
Letting $R\to\infty$ in \eqref{eq:observable-cutoff-identity} therefore yields
\(\Phi_\phi(t)-\Phi_\phi(s)
=\int_s^t\int L_{m(\tau)}\phi\dd\rho_\tau\dd\tau\).
The integrand is continuous in $\tau$, by narrow continuity of $\rho$
and continuity of $\tau\mapsto L_{m(\tau)}\phi$ in the uniform norm.
The fundamental theorem of calculus gives \eqref{eq:observable},
including the one-sided derivatives at $0$ and $T$.
\end{proof}

\begin{proposition}[Regularity and the consensus bootstrap]
\label[proposition]{prop:observables}
The following assertions are valid without regularity of the
initial measure.
\begin{enumerate}[label=\textup{(\roman*)},leftmargin=*]
\item For a prescribed $m\in C^r(I)$, with integer $r\geq0$, every fixed
$\phi\in\mathcal X^{2r+2}$ satisfies
$\Phi_\phi\in C^{r+1}(I)$. 
\item Suppose $w>0$ everywhere and
$w,v_1w,\ldots,v_dw\in\mathcal X^{2k}$, $k\geq1$.
Then $m\in C^k([0,T];\R^d)$.
\end{enumerate}
\end{proposition}

\begin{proof}
We prove (i) by induction on $r$.
For $r=0$, the conclusion is \cref{lem:observable}.
Let $r\geq1$ and assume that (i) has been proved at index $r-1$:
for each prescribed trajectory in $C^{r-1}(I)$, the integral against
$\rho_t$ of any fixed function in $\mathcal X^{2r}$ belongs to $C^r(I)$.
We now take $m\in C^r(I)$ and $\phi\in\mathcal X^{2r+2}$ and prove
that $\Phi_\phi\in C^{r+1}(I)$.

Apply \cref{lem:observable} to the fixed test function $\phi$ to obtain
$\Phi_\phi'(t)=\int L_{m(t)}\phi\dd\rho_t$.
Expanding the generator gives
\begin{equation}\label{eq:expanded-generator}
\begin{split}
L_m\phi={}&\frac{\sigma^2}{2}|v|^2\Delta\phi
                 -\lambda v\cdot\nabla\phi-\sigma^2\sum_{j=1}^dm_jv_j\Delta\phi
       +\frac{\sigma^2}{2}|m|^2\Delta\phi
       +\lambda\sum_{j=1}^dm_j\partial_j\phi.
\end{split}
\end{equation}
The spatial functions $|v|^2\Delta\phi$, $v\cdot\nabla\phi$,
$v_j\Delta\phi$, $\Delta\phi$, and $\partial_j\phi$ in
\eqref{eq:expanded-generator} all belong to $\mathcal X^{2r}$.
Each is a fixed function of $v$, independent of $t$ and $m$.
For example, $\Phi_{v_j\Delta\phi}(t)=\int v_j\Delta\phi(v)\dd\rho_t(v)$;
$m_j(t)$ multiplies this observable and is not part of its spatial test function.
Since $m\in C^r(I)\subset C^{r-1}(I)$, the induction hypothesis
applies to each of these functions: its integral against $\rho_t$
belongs to $C^r(I)$.

Substituting \eqref{eq:expanded-generator} into \eqref{eq:observable} gives
\begin{equation}\label{eq:observable-expanded}
\Phi_\phi'
=\frac{\sigma^2}{2}\Phi_{|v|^2\Delta\phi}
 -\lambda\Phi_{v\cdot\nabla\phi}
 -\sigma^2\sum_{j=1}^d m_j\Phi_{v_j\Delta\phi}
 +\frac{\sigma^2}{2}|m|^2\Phi_{\Delta\phi}
 +\lambda\sum_{j=1}^d m_j\Phi_{\partial_j\phi}.
\end{equation}
Every observable on the right side is $C^r(I)$ by the induction hypothesis,
and the coefficients $m_j(t)$ and $|m(t)|^2$ are $C^r(I)$ by assumption.
Thus every term on the right side of \eqref{eq:observable-expanded}
is a product of $C^r(I)$ functions, and the sum is finite.
It follows that $\Phi_\phi'\in C^r(I)$, hence
$\Phi_\phi\in C^{r+1}(I)$. This proves (i) at index $r$.

For (ii), we first establish continuity of $m$.
The weights $w$ and $v_iw$ are bounded and continuous, so narrow
continuity of $\rho_t$ implies that $B(t)=\Phi_w(t)$ and each
$A_i(t)=\Phi_{v_iw}(t)$ are continuous.
Since $w>0$ and $\rho_t$ is a probability measure, $B(t)>0$ for every
$t\in[0,T]$. Continuity on this compact interval gives
$\min_{[0,T]}B>0$. Hence each $m_i=A_i/B$ is continuous, and
$m\in C^0([0,T];\R^d)$. Because $k\geq1$, the weights $w,v_iw\in\mathcal X^{2k}$ belong
to $\mathcal X^2$. Apply (i) with $r=0$ to the continuous trajectory
$m$, first with the fixed test function $\phi=w$ and then with
$\phi=v_iw$. This gives $B=\Phi_w\in C^1([0,T])$ and
$A_i=\Phi_{v_iw}\in C^1([0,T])$ for every $i$.
The quotient rule gives $m_i'=(A_i'B-A_iB')/B^2$, which is continuous
because $B$ is bounded away from zero. Thus $m\in C^1([0,T];\R^d)$. Now let $1\leq j<k$ and assume that $m\in C^j([0,T];\R^d)$ has
already been proved. Since $2j+2\leq2k$, the fixed weights $w$ and
$v_iw$ belong to $\mathcal X^{2j+2}$.
Apply (i) with $r=j$ to these same test functions. We obtain
$B=\Phi_w\in C^{j+1}([0,T])$ and
$A_i=\Phi_{v_iw}\in C^{j+1}([0,T])$ for every $i$.
Since $\min_{[0,T]}B>0$, the reciprocal $1/B$ belongs to
$C^{j+1}([0,T])$. Consequently $m_i=A_i(1/B)\in C^{j+1}([0,T])$
for every $i$, so $m\in C^{j+1}([0,T];\R^d)$.
Starting with the $C^1$ regularity proved above and repeating this step
yields $m\in C^k([0,T];\R^d)$.
\end{proof}

\subsection{Local regularity away from the consensus point}
\label{sec:off-center}

For an open set $U$, write $C^{\infty,r}(U)$ for functions
whose mixed derivatives $\partial_v^\beta\partial_t^j$ are
continuous for every finite $\beta$ and $0\leq j\leq r$.
The following regularity result away from the consensus point is a
standard consequence of Schauder theory.

\begin{lemma}[Regularity away from the center]
\label[lemma]{lem:off-center}
Suppose $m\in C^r(I)$ for an integer $r\geq1$.
On the open set
$\{(v,t)\in\R^d\times I:v\neq m(t)\}$,
the space--time measure has a density $\rho\in C^{\infty,r+1}$.
\end{lemma}

\begin{proof}
Fix $(v_0,t_0)\in\R^d\times I$ with $v_0\neq m(t_0)$,
and set $R=|v_0-m(t_0)|/4$.
By continuity of $m$, choose $h>0$ such that
$[t_0-h,t_0+h]\Subset I$ and
$\sup_{|t-t_0|\leq h}|m(t)-m(t_0)|\leq R$.
The cylinder used below is
\(Q=B_R(v_0)\times(t_0-h,t_0+h)\).
For every $(v,t)\in\overline Q$, the triangle inequality gives
\(2R\leq |v-m(t)|\leq6R\).
Thus $Q\Subset\{(v,t)\in\R^d\times I:v\neq m(t)\}$.
In FPK notation put
\(A^{ij}(v,t)=\frac{\sigma^2}{2}|v-m(t)|^2\delta_{ij},  b^i(v,t)=-\lambda(v_i-m_i(t))\).
Then
$\partial_t\rho=\partial_{ij}(A^{ij}\rho)-\partial_i(b^i\rho)$.
On $Q$, $A,A^{-1},b$ are bounded and
\(\partial_{v_\ell}A^{ij} =\sigma^2(v_\ell-m_\ell(t))\delta_{ij}\).
Consequently, for every $1\leq p<\infty$,
\(\sup_{|t-t_0|\leq h}\|A(\cdot,t)\|_{W^{1,p}(B_R(v_0))} \leq C_{Q,p}<\infty\).
This $C_{Q,p}$ depends only on $d,p,\sigma,R$, since $|v-m(t)|\leq6R$ on $Q$.
Here the bound is local in both space and time: it is uniform
only for $t\in[t_0-h,t_0+h]$, and its constant may
depend on $Q$. No global $W^{1,p}(\R^d)$ bound or bound
uniform over all $t>0$ is asserted.
In particular we may choose $p>d+2$, and $b\in L^p(Q)$.
The coefficients are locally integrable against
$\rho_t\dd t$ by boundedness and the mass-one property.

Apply the parabolic measure regularity result
\cite[Corollary~6.4.3]{BKRS}.
It gives a density $\rho$ locally H\"older continuous in both space
and time, together with local first-order spatial Sobolev
regularity; in particular, $\rho,\nabla_v\rho\in L^p_{\rm loc}(Q)$.
The time H\"older continuity follows explicitly from the
parabolic embedding step in the proof of
\cite[Theorem~6.4.2]{BKRS}.
The preliminary density existence
result is \cite[Corollary~6.3.2]{BKRS}.
These results start with a measure solution, so no function
density has been assumed at this step.
For any $\psi\in C_c^\infty(B_R(v_0))$, the identity
$\int\psi\dd\rho_t=\int\psi(v)\rho(v,t)\dd v$ holds for almost every
$t\in(t_0-h,t_0+h)$. Both sides are continuous in $t$, by narrow
continuity on the left and continuity of the density on the right,
so the identity holds at every time. Thus $\rho$ represents every
marginal measure locally. The densities on overlapping cylinders agree
everywhere by continuity, since they agree almost everywhere.
On each smaller cylinder, the resulting density satisfies
\(|\rho(v,t)-\rho(v',t')| \leq C\big(|v-v'|^\theta+|t-t'|^{\theta/2}\big)\)
for all pairs $(v,t),(v',t')$ in that cylinder and some
$\theta\in(0,1)$, decreasing the exponent if
necessary. Here $C$ and $\theta$ depend on $d,p,\lambda,\sigma,Q$, the smaller cylinder and the mass of $\rho_t\dd t$ on $Q$.
This is local parabolic H\"older continuity,
obtained from the measure equation; it is not a regularity
assumption on the initial measure $\rho_0$.
The spatial Sobolev assertion concerns weak derivatives
in $v$; it is distinct from this space--time H\"older estimate.

The divergence form, including its first-order coefficient, is
\begin{equation}\label{eq:divergence}
\partial_t\rho
 =\divv_v\left(\frac{\sigma^2}{2}|v-m(t)|^2\nabla_v\rho
                  +(\lambda+\sigma^2)(v-m(t))\rho\right).
\end{equation}
Indeed, with $y=v-m(t)$,
$\Delta(|y|^2\rho)=\divv(|y|^2\nabla \rho+2y\rho)$.
The forward operator in \eqref{eq:forward} is
\(L_m^*=\frac{\sigma^2}{2}|v-m|^2\Delta +(\lambda+2\sigma^2)(v-m)\cdot\nabla +d(\lambda+\sigma^2)\).
Since $m\in C^1$ (by \cref{prop:observables}(ii) for the nonlinear
equation) and the coefficients are polynomial in $v$, they
and all their spatial derivatives are locally
$C^{\theta,\theta/2}$.

We justify the passage from a weak solution to a classical
one before applying a classical Schauder estimate.
Here the local energy conditions mean that, for every spatial
ball $B$ and interval $[a,b]$ with
$\overline B\times[a,b]\Subset Q$,
\(\sup_{t\in[a,b]}\int_B|\rho(v,t)|^2\dd v +\int_a^b\int_B|\nabla_v\rho(v,t)|^2\dd v\dd t<\infty\).
The first bound, and $\rho\in C([a,b];L^2(B))$, follow from
the continuity of the density on the compact cylinder.
The second follows from $\nabla_v\rho\in L^p_{\rm loc}(Q)$
with $p>d+2>2$. Thus these are consequences of the preceding
measure regularity result, not additional assumptions.
Apply \cite[Theorem~6.6]{Lieberman} to \eqref{eq:divergence}.
All spatial derivatives of the coefficients are locally bounded,
so that result and its stated induction give
\(\rho,\partial_t\rho\in L^2((a,b);H^N(B)) \quad\text{for every integer }N\geq0\),
on every such cylinder. Thus
$\rho\in H^1((a,b);H^N(B))\subset C^{1/2}([a,b];H^N(B))$.
Taking $N$ arbitrarily large and using spatial Sobolev
embedding shows that every spatial derivative is jointly
continuous and locally $1/2$-H\"older in time.
Since $m\in C^1$, the distributional identity
$\partial_t\rho=L_{m(t)}^*\rho$ now has a continuous right-hand
side, also after every spatial derivative. Integration in
time makes these time derivatives classical.
In particular, $\nabla_v^2\rho$ and $\partial_t\rho$ are locally
$C^{\theta,\theta/2}$. The parabolic norm characterization in
\cite[Corollary~4.10]{Lieberman}
then gives $\rho\in C^{2+\theta,1+\theta/2}$ locally.

The precise classical interior Schauder estimate is
\cite[Theorem~4.9]{Lieberman}.
It now applies to $\rho$ and its spatial derivatives.
Their equations are
\begin{equation}\label{eq:schauder-recursion}
\begin{split}
(\partial_t-L_{m(t)}^*)
        \partial_v^\beta \rho
 =\sum_{0<\gamma\leq\beta}\binom{\beta}{\gamma}
 \big[&\tfrac{\sigma^2}{2}(\partial_v^\gamma|y|^2)
               \Delta\partial_v^{\beta-\gamma}\rho\\
     &+(\lambda+2\sigma^2)(\partial_v^\gamma y)\cdot
                \nabla\partial_v^{\beta-\gamma}\rho\big].
\end{split}
\end{equation}
The weak-to-classical step above is supplied by
\cite[Theorem~6.6]{Lieberman};
\cite[Theorem~4.9]{Lieberman} gives interior H\"older
bounds once that step has been established.
Repeating the time differentiation argument, also after any spatial derivative,
gives
\begin{equation}\label{eq:time-bootstrap}
\partial_t^{j+1}\rho
 =\sum_{\ell=0}^{j}\binom j\ell
       (\partial_t^\ell L_{m(t)}^*)\partial_t^{j-\ell}\rho,
\qquad 0\leq j\leq r.
\end{equation}
Only derivatives of $m$ through order $r$ are used.
Since $(v_0,t_0)$ was arbitrary, this proves
$C^{\infty,r+1}$ regularity on the stated open set.
\end{proof}

\section{Hypoelliptic estimates and the proof of Theorem~\ref{thm:main}}
\label{sec:finite-full}

We begin with $m\in C^1$ from \cref{prop:observables}(ii) and first
treat $m$ as a prescribed trajectory. In \cref{sec:fourier-equation},
we prove a fourth-moment bound for the given weak solution and derive
the equation for the centered Fourier transform $F$.
The dilation $F_R(t,\xi)=F(t,R\xi)$ then puts every high-frequency
annulus on a fixed annulus, with the scale $R$ appearing in the
imaginary potential $-iR\,m'(t)\cdot\xi$. Nonconstant motion provides an interval on which one component of $m'$
is bounded away from zero. On this interval,
\cref{sec:annular-commutator} proves an $R^{2/3}$ gain by combining
Fourier energy estimates with a commutator in frequency.
The constants are uniform in $R$, and the argument requires only
continuity of $m'$, since no time derivative of $m'$ is taken. In \cref{sec:fourier-localization}, we localize this gain to control
the solution on an inner cylinder by $R^{-2/3}$ times its norm on a
larger cylinder. Starting from the probability bound $|F_R|\leq1$,
iteration along nested cylinders gives arbitrarily high powers of
$R^{-1}$ in the annular bounds. Dyadic summation then yields every
spatial $H^s$ order on the interval of motion. In \cref{sec:forward-propagation}, a weighted Fourier energy estimate
propagates each acquired $H^s$ bound to every later time, including
times when $m$ has stopped moving. Every time after $t_c$ has an
earlier interval of nonconstant motion, so this gives all spatial
Sobolev orders throughout $(t_c,T)$. In \cref{sec:post-motion-bootstrap}, we return to $m=A/B$.
Once all spatial Sobolev orders are available, the PDE gives
$\rho\in C^{\infty,r+1}$ whenever $m\in C^r$, $r\geq1$.
Rapid decay of the fixed weights justifies differentiating the
integrals defining $A$ and $B$ in time, giving $A,B\in C^{r+1}$ and
hence $m\in C^{r+1}$. Iteration starting at $r=1$ makes $m$ smooth
in time and $\rho$ jointly smooth in space and time after $t_c$. Finally, \cref{sec:onset-completion} combines these conclusions with
the local regularity away from the consensus point in
\cref{lem:off-center}. This treats the initial stationary phase and
an interior transition time, and completes the regularity proof of
\cref{thm:main}.

\subsection{Fourth moments and the exact Fourier equation}
\label{sec:fourier-equation}

The proofs of \cref{lem:fourth-moment,lem:fourier-equation} are given in \cref{sec:moment-fourier-proofs}.

\begin{lemma}[Moment propagation for the given weak solution]
\label[lemma]{lem:fourth-moment}
Suppose $m\in C([0,T];\R^d)$ and $\rho_0\in\Pp_4$.
Every narrowly continuous probability solution of the linear
equation with this prescribed $m$ satisfies
\begin{equation}\label{eq:fourth-moment}
\sup_{0\leq t\leq T}\int(1+|v|^2)^2\dd\rho_t(v)
 \leq e^{CT}\int(1+|v|^2)^2\dd\rho_0(v),
\end{equation}
where $C$ depends only on $d,\lambda,\sigma$ and $\sup_{[0,T]}|m|$. No moment bound on the given solution
is presupposed.
\end{lemma}

\begin{lemma}[Fourier equation for a measure solution]
\label[lemma]{lem:fourier-equation}
Let $m\in C^1([0,T];\R^d)$ and $\rho_0\in\Pp_4(\R^d)$.
Let $\rho$ be a narrowly continuous probability-measure
solution of \eqref{eq:weak} with this prescribed $m$
and initial value $\rho_0$. Define
\begin{equation}\label{eq:moving-fourier}
F(t,\xi)
 =\int e^{-i\xi\cdot(v-m(t))}\dd\rho_t(v)
 =e^{i\xi\cdot m(t)}\widehat\rho_t(\xi),
\qquad |F|\leq1.
\end{equation}
Then $F$, its first time derivative, and its frequency
derivatives through order two are jointly continuous on
$[0,T]\times\R^d$, with one-sided time derivatives at
the endpoints. Moreover,
\begin{equation}\label{eq:fourier-equation}
\partial_tF
 =\frac{\sigma^2}{2}|\xi|^2\Delta_\xi F
       -\lambda\xi\cdot\nabla_\xi F
       +i\,m'(t)\cdot\xi\,F.
\end{equation}
\end{lemma}

The function $F$ is the Fourier transform in coordinates centered at
the consensus point. The phase factor shifts $v$ to
$v-m(t)$ and cancels the terms containing $m(t)$ itself
in the equation for $\widehat\rho_t$, leaving only the velocity $m'(t)$
in \eqref{eq:fourier-equation}. Since the phase has modulus one,
$|F|=|\widehat\rho_t|$, so weighted Fourier $L^2$ estimates for $F$
give exactly the same spatial Sobolev norms of $\rho_t$. This form has two useful features. First, the trajectory
enters only through its velocity: after frequency dilation
by $R$, the potential becomes $iR\,m'(t)\cdot\xi$.
Differentiating it in $\xi_j$ produces $iR m_j'(t)$,
which is the large term in the commutator argument
leading to the $R^{2/3}$ gain in
\eqref{eq:annular-subelliptic}. Second, this potential
is purely imaginary. Multiplication by a real cutoff
and pairing with the complex conjugate of the solution
therefore give zero contribution from it to the real-part
energy identity. This is why the local energy estimate
\eqref{eq:annular-caccioppoli} has a constant independent
of $R$ and $m'$. These arguments differentiate the
potential only in frequency and require no $m''$.

For $R\geq1$, define the frequency dilation
$F_R(t,\xi)=F(t,R\xi)$. Then
\begin{equation}\label{eq:scaled-fourier}
\mathcal L_R F_R=0,\qquad
\mathcal L_R
 =\partial_t-\frac{\sigma^2}{2}|\xi|^2\Delta_\xi 
          +\lambda \xi\cdot\nabla_\xi -iR\,m'(t)\cdot \xi.
\end{equation}
On a fixed annulus the diffusion is uniformly elliptic,
while the motion of the consensus appears as an imaginary
linear potential whose size is $R$.

\subsection{The commutator estimate, with no derivative of the velocity}
\label{sec:annular-commutator}

\begin{lemma}[Uniform annular subelliptic estimate]
\label[lemma]{lem:annular-subelliptic}
Let $I$ be a bounded open interval and
$A=\{\xi:r_-<|\xi|<r_+\}$, where $0<r_-<r_+<\infty$.
Suppose $m\in C^1(I;\R^d)$ with bounded $m'$, and for some
$j\in\{1,\ldots,d\}$ and $\delta>0$,
\begin{equation}\label{eq:velocity-sign}
|m_j'(t)|\geq\delta\qquad(t\in I).
\end{equation}
By continuity, $m_j'$ has a fixed sign on $I$.
There is a constant $C$, independent of $R\geq1$, such that
\begin{equation}\label{eq:annular-subelliptic}
R^{2/3}\|U\|_{L^2(I\times A)}
 \leq C\big(\|\mathcal L_RU\|_{L^2(I\times A)}
                         +\|U\|_{L^2(I\times A)}\big)
\end{equation}
for every complex-valued function $U$ compactly supported in
$I\times A$ such that $U$, $\partial_tU$, $\nabla_\xi U$
and $\nabla_\xi^2U$ are continuous.
The constant depends only on $r_-,r_+,d,\lambda,\sigma,\|m'\|_{L^\infty(I)}$ and $\delta$; in particular it does not depend on the length of $I$.
\end{lemma}

\begin{proof}
We first prove the estimate for $U\in C_c^\infty(I\times A)$,
then pass to the stated function class by mollification.
Every $C$ and the positive $c$ in this proof depend only
on $d,\lambda,\sigma,r_-,r_+,\|m'\|_{L^\infty(I)}$
and $\delta$. No derivative or modulus of continuity of $m'$ is used.
All norms and inner products in this proof are over
$I\times\R^d$, with $U$ extended by zero. The inner product
is linear in its first argument. We write $\nabla_\xi^2U$
for the Hessian matrix. The identity
$\|\Delta_\xi U\|_2=\|\nabla_\xi^2U\|_2$ follows from
the Fourier transform in $\xi$.

First, take the real part of the pairing with $U$.
Compact support in time gives
$\operatorname{Re}\int(\partial_tU)\overline U
=\frac12\int\partial_t|U|^2=0$.
The potential contributes zero because
$-iR(m'(t)\cdot\xi)|U|^2$ is purely imaginary.
For the diffusion, use
$\operatorname{Re}(\overline U\,\nabla_\xi U)
=\frac12\nabla_\xi|U|^2$ and integrate by parts:
\[
\begin{aligned}
-\frac{\sigma^2}{2}\operatorname{Re}
  \int|\xi|^2(\Delta_\xi U)\overline U
 &=\frac{\sigma^2}{2}\int|\xi|^2|\nabla_\xi U|^2
   +\sigma^2\operatorname{Re}
       \int\xi\cdot\nabla_\xi U\,\overline U\\
 &=\frac{\sigma^2}{2}\int|\xi|^2|\nabla_\xi U|^2
   +\frac{\sigma^2}{2}\int\xi\cdot\nabla_\xi|U|^2\\
 &=\frac{\sigma^2}{2}\int|\xi|^2|\nabla_\xi U|^2
   -\frac{d\sigma^2}{2}\|U\|_2^2.
\end{aligned}
\]
Here $\nabla_\xi|\xi|^2=2\xi$ and
$\operatorname{div}_\xi\xi=d$; all spatial boundary terms
vanish by compact support. The drift similarly gives
\[
\lambda\operatorname{Re}\int
 (\xi\cdot\nabla_\xi U)\overline U
 =\frac\lambda2\int\xi\cdot\nabla_\xi|U|^2
 =-\frac{d\lambda}{2}\|U\|_2^2.
\]
Adding these contributions yields
\(\operatorname{Re}\langle\mathcal L_RU,U\rangle =\frac{\sigma^2}{2}\int|\xi|^2|\nabla_\xi U|^2 -\frac{d(\lambda+\sigma^2)}2\|U\|_2^2\).
Thus, using $r_-\leq |\xi|\leq r_+$ on the support,
\begin{equation}\label{eq:annular-G}
\begin{aligned}
\|\nabla_\xi U\|_2^2
 &\leq C\big(\|\mathcal L_RU\|_2+\|U\|_2\big)\|U\|_2 ,\\
\|\nabla_\xi U\|_2
 &\leq C\big(\|\mathcal L_RU\|_2+\|U\|_2\big).
\end{aligned}
\end{equation}

Second, pair the equation with $-\Delta_\xi U$ and take real
parts. The time derivative integrates to zero. The
second-order part contributes exactly
$(\sigma^2/2)\int|\xi|^2|\Delta_\xi U|^2$.
For the drift, the product rule gives
$\partial_{\xi_j}(\xi\cdot\nabla_\xi U)
=\partial_{\xi_j}U+\xi\cdot\nabla_\xi\partial_{\xi_j}U$.
Hence integration by parts yields
\[
\begin{aligned}
\operatorname{Re}\langle\lambda\xi\cdot\nabla_\xi U,
                         -\Delta_\xi U\rangle
 &=\lambda\sum_j\operatorname{Re}\int
   \partial_{\xi_j}(\xi\cdot\nabla_\xi U)
          \overline{\partial_{\xi_j}U}\\
 &=\lambda\|\nabla_\xi U\|_2^2
   +\lambda\sum_j\operatorname{Re}\int
    (\xi\cdot\nabla_\xi\partial_{\xi_j}U)
          \overline{\partial_{\xi_j}U}\\
 &=\lambda\|\nabla_\xi U\|_2^2
   +\frac\lambda2\int\xi\cdot\nabla_\xi|\nabla_\xi U|^2\\
 &=\lambda(1-d/2)\|\nabla_\xi U\|_2^2.
\end{aligned}
\]
The last equality again uses $\operatorname{div}_\xi\xi=d$
and compact support. For the potential, $m'(t)$ is independent
of $\xi$, so a single integration by parts gives
\[
\begin{aligned}
\langle-iR(m'(t)\cdot\xi)U,-\Delta_\xi U\rangle
 &=iR\int(m'(t)\cdot\xi)U\,\Delta_\xi\overline U\\
 &=-iR\int U\,m'(t)\cdot\nabla_\xi\overline U
   -iR\int(m'(t)\cdot\xi)|\nabla_\xi U|^2.
\end{aligned}
\]
The second term is purely imaginary. Taking real parts
and applying Cauchy--Schwarz to the first term therefore yields
\[
\begin{aligned}
\left|\operatorname{Re}
 \langle-iR(m'(t)\cdot\xi)U,-\Delta_\xi U\rangle\right|
 &\leq R\int|m'(t)|\,|U|\,|\nabla_\xi U|\\
 &\leq R\|m'\|_{L^\infty(I)}
              \|U\|_2\|\nabla_\xi U\|_2.
\end{aligned}
\]
Consequently,
\[
\begin{aligned}
\|\Delta_\xi U\|_2^2
 &\leq C\big(\|\mathcal L_RU\|_2+\|U\|_2\big)
              \|\Delta_\xi U\|_2
       +C\|\nabla_\xi U\|_2^2
       +CR\|U\|_2\|\nabla_\xi U\|_2.
\end{aligned}
\]
Using \eqref{eq:annular-G} and Young's inequality gives
\[
\|\Delta_\xi U\|_2
 \leq C\big(\|\mathcal L_RU\|_2+\|U\|_2
       +R^{1/2}\|U\|_2^{1/2}\|\nabla_\xi U\|_2^{1/2}\big).
\]

Third, introduce the skew-adjoint operator
\(\mathcal T_R=\partial_t+\lambda \xi\cdot\nabla_\xi +\frac{d\lambda}{2}-iR\,m'(t)\cdot \xi,  \mathcal T_R^*=-\mathcal T_R\).
The frequency derivative $\partial_{\xi_j}$ is also
skew-adjoint, and the exact commutator is
\begin{equation}\label{eq:fourier-bracket}
[\partial_{\xi_j},\mathcal T_R]
 =\lambda\partial_{\xi_j}-iR m_j'(t).
\end{equation}
In particular, this computation takes no time derivative
of $m'$. Skew-adjointness gives
\(\langle[\partial_{\xi_j},\mathcal T_R]U,U\rangle =\langle\partial_{\xi_j}U,\mathcal T_RU\rangle -\langle\mathcal T_RU,\partial_{\xi_j}U\rangle\),
whose absolute value is at most
$2\|\nabla_\xi U\|_2\|\mathcal T_RU\|_2$.
On the other hand, \eqref{eq:velocity-sign} and
\eqref{eq:fourier-bracket} imply
\[
R\delta\|U\|_2^2
 \leq2\|\nabla_\xi U\|_2\|\mathcal T_RU\|_2
       +\lambda\|\nabla_\xi U\|_2\|U\|_2.
\]
As
\(\mathcal T_RU=\mathcal L_RU +\frac{\sigma^2}{2}|\xi|^2\Delta_\xi U+\frac{d\lambda}{2}U\),
we obtain
\[
\begin{aligned}
R\|U\|_2^2
 &\leq C\|\nabla_\xi U\|_2
       \big(\|\mathcal L_RU\|_2+\|U\|_2+\|\Delta_\xi U\|_2\big)\\
 &\leq C\big[\big(\|\mathcal L_RU\|_2+\|U\|_2\big)
                     \|\nabla_\xi U\|_2
       +R^{1/2}\|U\|_2^{1/2}\|\nabla_\xi U\|_2^{3/2}\big].
\end{aligned}
\]
From here on $C$ also depends on $\delta$, through the division of the commutator bound by the velocity lower bound.
Substituting \eqref{eq:annular-G} yields
\[
R\|U\|_2^2
 \leq C\big[\big(\|\mathcal L_RU\|_2+\|U\|_2\big)^{3/2}
                     \|U\|_2^{1/2}
       +R^{1/2}\big(\|\mathcal L_RU\|_2+\|U\|_2\big)^{3/4}
                     \|U\|_2^{5/4}\big].
\]
If $\|U\|_2>0$, division by $R\|U\|_2^2$ gives
\[
1\leq C\left[
 \left(\frac{\|\mathcal L_RU\|_2+\|U\|_2}
                  {R^{2/3}\|U\|_2}\right)^{3/2}
 +\left(\frac{\|\mathcal L_RU\|_2+\|U\|_2}
                  {R^{2/3}\|U\|_2}\right)^{3/4}\right].
\]
Thus $\|\mathcal L_RU\|_2+\|U\|_2
\geq cR^{2/3}\|U\|_2$ for a constant $c>0$.
This is \eqref{eq:annular-subelliptic}; $\|U\|_2=0$ is immediate.

Finally, let $U$ belong to the function class in the statement.
Extend it by zero and mollify in $(t,\xi)$.
Since its support is a compact subset of $I\times A$,
sufficiently small smoothing radii give
$U_n\in C_c^\infty(I\times A)$ with supports in a common
compact subset. The functions and their first time and
first two spatial derivatives converge uniformly.
For each fixed $R$, this implies $U_n\to U$ and
$\mathcal L_RU_n\to\mathcal L_RU$ in $L^2$.
Applying the estimate to $U_n$ and passing to the limit
proves \eqref{eq:annular-subelliptic} for $U$ with the same constant. Only $U$ is mollified;
$m'$ remains unchanged and is never differentiated in time.
\end{proof}

Related commutator estimates were established by Bouchut
\cite[Theorems~1.5--1.6 and Section~3]{Bouchut2002}.
For the kinetic Fokker--Planck equation
$\partial_t\phi+v\cdot\nabla_x\phi-\Delta_v\phi=h$,
with $\phi,h,\nabla_v\phi\in L^2$, he established
$\|\partial_t\phi+v\cdot\nabla_x\phi\|_2
+\|\Delta_v\phi\|_2+\|D_x^{2/3}\phi\|_2
\leq C\|h\|_2$.
His proof uses the identity
$[\partial_{v_j},\partial_t+v\cdot\nabla_x]=\partial_{x_j}$
to relate spatial derivatives to velocity derivatives and the
transport operator through integration by parts.
Energy estimates control the diffusive derivatives, and Fourier
interpolation closes the resulting inequalities.
This combination of commutator identities, energy estimates,
and interpolation is also the mechanism underlying our annular
estimate, although the operator and the quantities being
controlled are different.

Our contribution is to establish a commutator estimate uniform
in the frequency scale for the CBO equation, using only the
time regularity available before the nonlinear bootstrap.
The kinetic estimates above do not directly give the required
bound for the scaled operator $\mathcal L_R$, whose diffusion
coefficient is $|\xi|^2$ and whose imaginary potential contains
the continuous velocity $m'(t)$.
The identity
$[\partial_{\xi_j},\mathcal T_R]
=\lambda\partial_{\xi_j}-iRm_j'(t)$,
together with annular energy estimates, yields
$R^{2/3}\|U\|_2
\leq C(\|\mathcal L_RU\|_2+\|U\|_2)$
with $C$ independent of $R$, without differentiating $m'$ in
time (\cref{lem:annular-subelliptic}).

A further step~(see section~\ref{sec:fourier-localization}) is to localize this estimate uniformly in $R$
and iterate it on nested cylinders.
The probability bound $|F_R|\leq1$ then supplies the starting
control, allowing the iteration to produce every spatial
$H^s$ order from measure initial data, without an initial
density or $L^2$ assumption (\cref{prop:fourier-smoothing}).

\subsection{Localization and iteration to every Sobolev order}
\label{sec:fourier-localization}

\begin{lemma}[Localized annular gain]
\label[lemma]{lem:localized-annular-gain}
Let $I,A,m$ and $\mathcal L_R$ be as in
\cref{lem:annular-subelliptic}, and set $Q=I\times A$.
Fix product cylinders $Q_0\Subset Q_1\Subset Q_2\Subset Q$,
with $Q_0=I_0\times A_0$.
For every $R\geq1$, let $U$ solve $\mathcal L_RU=0$
on $Q_2$, with $U,\partial_tU,\nabla_\xi U,\nabla_\xi^2U$
continuous there. Then
\begin{equation}\label{eq:one-annular-gain}
\begin{aligned}
\|U\|_{L^2(Q_0)}
 &\leq C R^{-2/3}\|U\|_{L^2(Q_2)},\\
\sup_{t\in I_0}\int_{A_0}|U(t,\xi)|^2\dd\xi
 &\leq C R^{-4/3}\|U\|_{L^2(Q_2)}^2.
\end{aligned}
\end{equation}
Here $C$ depends only on $d,\lambda,\sigma,r_-,r_+,\delta,\|m'\|_{L^\infty(I)}$ and the geometry of $Q_0\Subset Q_1\Subset Q_2$; it does not depend on $R$.
\end{lemma}

\begin{proof}
We first prove the local energy estimate used below.
Write $Q_2=I_2\times A_2$. For every product cylinder
$Q'=I'\times A'\Subset Q_2$, we claim
\begin{equation}\label{eq:annular-caccioppoli}
\sup_{t\in I'}\int_{A'}|U(t,\xi)|^2\dd \xi
       +\int_{Q'}|\nabla_\xi U|^2\dd \xi\dd t
 \leq C_{Q',Q_2}\int_{Q_2}|U|^2\dd \xi\dd t,
\end{equation}
where $C_{Q',Q_2}$ depends only on $d,\lambda,\sigma,r_-,r_+$ and the $C^1$ norms of the cutoffs chosen below; it does not depend on $R$ or on $m'$.
Let $a$ be the left endpoint of $I_2$ and
choose real cutoffs $\chi\in C_c^\infty(I_2)$ and
$\theta\in C_c^\infty(A_2)$, taking values in $[0,1]$,
with $\chi=1$ on $I'$ and $\theta=1$ on $A'$.
At each fixed time, multiply $\mathcal L_RU=0$ by
$\chi^2\theta^2\overline U$, integrate in $\xi$, and take
real parts. The potential contributes zero because
$-iR(m'(t)\cdot\xi)\chi^2\theta^2|U|^2$ is purely imaginary.
For the diffusion use
$\nabla_\xi(\theta^2|\xi|^2)
=2\theta|\xi|^2\nabla_\xi\theta+2\theta^2\xi$;
for the drift use
$\operatorname{div}_\xi(\theta^2\xi)
=d\theta^2+2\theta\xi\cdot\nabla_\xi\theta$.
Integration by parts therefore gives the exact identity
\[
\begin{aligned}
&\frac12\frac{\dd}{\dd t}
 \left(\chi^2\int_{A_2}\theta^2|U|^2\dd\xi\right)
 +\frac{\sigma^2}{2}\chi^2
       \int_{A_2}\theta^2|\xi|^2|\nabla_\xi U|^2\dd\xi\\
 &\qquad=\chi\chi'\int_{A_2}\theta^2|U|^2\dd\xi\\
 &\qquad\quad-\sigma^2\chi^2\operatorname{Re}\int_{A_2}
 (\theta|\xi|^2\nabla_\xi\theta+\theta^2\xi)
       \cdot\nabla_\xi U\,\overline U\dd\xi\\
 &\qquad\quad+\frac\lambda2\chi^2\int_{A_2}
 (d\theta^2+2\theta\xi\cdot\nabla_\xi\theta)|U|^2\dd\xi.
\end{aligned}
\]
Since $r_-<|\xi|<r_+$, Young's inequality bounds the
diffusion cross term pointwise by
\[
\begin{aligned}
\sigma^2\chi^2
 \left|(\theta|\xi|^2\nabla_\xi\theta+\theta^2\xi)
           \cdot\nabla_\xi U\,\overline U\right|
 &\leq\frac{\sigma^2}{4}\chi^2\theta^2|\xi|^2
                         |\nabla_\xi U|^2\\
 &\quad+C\chi^2(\theta^2+|\nabla_\xi\theta|^2)|U|^2.
\end{aligned}
\]
The drift term is bounded by
$C\chi^2(\theta^2+|\nabla_\xi\theta|^2)|U|^2$
under the integral, while the time-cutoff term is bounded
by $|\chi\chi'|\theta^2|U|^2$.
Absorb the first term of the preceding bound into the
left-hand side and integrate from $a$ to any $\tau\in I'$.
Since $\chi$ vanishes near $a$, the initial boundary term is
zero, and we obtain
\[
\begin{aligned}
\frac12\chi(\tau)^2\int_{A_2}\theta^2|U(\tau,\xi)|^2\dd\xi
 &+\frac{\sigma^2}{4}\int_a^\tau\int_{A_2}
       \chi^2\theta^2|\xi|^2|\nabla_\xi U|^2\dd\xi\dd t\\
 &\leq C_{Q',Q_2}\int_{Q_2}|U|^2\dd\xi\dd t.
\end{aligned}
\]
Taking the supremum over $\tau\in I'$ controls the first
term in \eqref{eq:annular-caccioppoli}. Letting $\tau$
increase to the right endpoint of $I'$ controls the
gradient integral over $Q'$, because $\chi=\theta=1$
there and $|\xi|\geq r_-$. Adding these two bounds proves
\eqref{eq:annular-caccioppoli}. Its constant is independent
of $R$ because the only term containing $R$ vanished
before any estimate was taken.

Now choose real cutoffs taking values in $[0,1]$, with
product $\chi(t)\theta(\xi)$ equal to one on $Q_0$
and supported in $Q_1$. Extend $\chi\theta U$ by zero
outside $Q_2$; its support lies strictly inside $Q_2$,
so this extension retains the required derivatives on $Q$.
There is no large potential term
in the cutoff commutator, since multiplication by
$-iR\,m'(t)\cdot\xi$ commutes with the cutoff.
Using $\mathcal L_RU=0$ and expanding the derivatives gives
\[
\begin{aligned}
\mathcal L_R(\chi\theta U)
 &=\chi'\theta U
   -\frac{\sigma^2}{2}\chi|\xi|^2
       \big((\Delta_\xi\theta)U
             +2\nabla_\xi\theta\cdot\nabla_\xi U\big)\\
 &\quad+\lambda\chi(\xi\cdot\nabla_\xi\theta)U.
\end{aligned}
\]
The coefficients on the right are bounded independently of $R$
and supported in $Q_1$. Thus
\[
\|\mathcal L_R(\chi\theta U)\|_{L^2(Q)}
 \leq C\big(\|U\|_{L^2(Q_1)}
                      +\|\nabla_\xi U\|_{L^2(Q_1)}\big).
\]
The function $\chi\theta U$ is compactly supported in $Q$
and has the regularity required by
\cref{lem:annular-subelliptic}. Since it equals $U$ on $Q_0$,
that lemma gives
\[
\begin{aligned}
R^{2/3}\|U\|_{L^2(Q_0)}
 &\leq R^{2/3}\|\chi\theta U\|_{L^2(Q)}\\
 &\leq C\big(\|\mathcal L_R(\chi\theta U)\|_{L^2(Q)}
                           +\|\chi\theta U\|_{L^2(Q)}\big)\\
 &\leq C\big(\|U\|_{L^2(Q_1)}
                      +\|\nabla_\xi U\|_{L^2(Q_1)}\big).
\end{aligned}
\]
Finally, \eqref{eq:annular-caccioppoli} on
$Q_1\Subset Q_2$ yields
$\|\nabla_\xi U\|_{L^2(Q_1)}\leq C\|U\|_{L^2(Q_2)}$,
and $\|U\|_{L^2(Q_1)}\leq\|U\|_{L^2(Q_2)}$ by inclusion.
Dividing by $R^{2/3}$ proves \eqref{eq:one-annular-gain} for the space--time norm.
To obtain the uniform-in-time estimate, first apply
\eqref{eq:annular-caccioppoli} on $Q_0\Subset Q_1$.
Next apply the space--time estimate just proved with
inner cylinder $Q_1$ and outer cylinder $Q_2$, choosing
an intermediate product cylinder between them. This gives
\[
\sup_{t\in I_0}\int_{A_0}|U(t,\xi)|^2\dd\xi
 \leq C\|U\|_{L^2(Q_1)}^2
 \leq C R^{-4/3}\|U\|_{L^2(Q_2)}^2.
\] This proves the second estimate
in \eqref{eq:one-annular-gain}.
Only multiplication by cutoffs and first frequency derivatives
of $U$ occurred. In particular the estimate can be iterated
without differentiating the equation in time.
\end{proof}

\begin{proposition}[Smoothing during a nonzero consensus velocity]
\label[proposition]{prop:fourier-smoothing}
Suppose \eqref{eq:velocity-sign} holds on an open interval
$I\Subset(0,T)$. For every $J\Subset I$ and every $s\geq0$,
\begin{equation}\label{eq:all-sobolev-good}
\sup_{t\in J}\|\rho_t\|_{H^s(\R^d)}<\infty,
\qquad \rho\in C(J;H^s(\R^d)).
\end{equation}
The conclusion holds for the given measure solution.
\end{proposition}

\begin{proof}
Apply \cref{lem:localized-annular-gain} to $F_R$.
For each fixed integer $n\geq1$ choose a finite chain
of nested cylinders between
$J\times\{1\leq|\xi|\leq2\}$ and
$I\times\{1/2<|\xi|<4\}$. Apply the space--time
estimate in \eqref{eq:one-annular-gain} $n-1$ times,
then its uniform-in-time estimate on the innermost
cylinders. Together with $|F_R|\leq1$, this gives
\(\sup_{t\in J}\int_{1\leq|\xi|\leq2}|F_R(t,\xi)|^2\dd \xi \leq C_{n,J}R^{-4n/3}, R\geq1\).
Here $C_{n,J}$ depends only on $d,\lambda,\sigma,n,\delta,\|m'\|_{L^\infty(I)}$ and the chosen cylinders and cutoffs; it does not depend on $R$, on $t\in J$ or on the initial measure.
Each of the first $n-1$ applications contributes
$R^{-4/3}$ to the squared space--time norm, and the
last contributes another $R^{-4/3}$ while taking the
time supremum. Their product is $R^{-4n/3}$, as displayed.
The last input is the explicit bound
\[
\int_I\int_{1/2<|\xi|<4}|F_R|^2\dd\xi\dd t
 \leq |I|\,|\{1/2<|\xi|<4\}|,
\]
which follows from $|F_R|\leq1$ and accounts for the
outer cylinder's volume in $C_{n,J}$.

The dilation relation gives
\[
\sup_{t\in J}\int_{R\leq|\xi|\leq2R}|F(t,\xi)|^2\dd\xi
 =R^d\sup_{t\in J}\int_{1\leq|\xi|\leq2}|F_R(t,\xi)|^2\dd\xi
       \leq C_{n,J}R^{d-4n/3}.
\]
Since $|F|=|\widehat\rho_t|$, this is also an annular
bound for $\widehat\rho_t$. Fix $s\geq0$ and choose an
integer $n>\frac34(2s+d)$; the preceding argument allows
any such $n$. Decompose frequency space into $|\xi|<1$
and the disjoint annuli
$2^\ell\leq|\xi|<2^{\ell+1}$, $\ell=0,1,\ldots$.
On the latter annulus,
$(1+|\xi|^2)^s\leq5^s2^{2s\ell}$.
Applying the preceding estimate with $R=2^\ell$ gives
\[
\sup_{t\in J}
 \int_{2^\ell\leq|\xi|<2^{\ell+1}}
       (1+|\xi|^2)^s|\widehat\rho_t(\xi)|^2\dd\xi
 \leq5^sC_{n,J}\,2^{\ell(2s+d-4n/3)}.
\]
The constant does not depend on $\ell$, because the annular estimate is uniform in $R\geq1$.
For $|\xi|<1$, the mass-one bound
$|\widehat\rho_t|=|F|\leq1$ gives an integral at most
$2^s|B_1|$. Summing the nonnegative integrals and using
that the supremum of a sum is bounded by the sum of the
suprema, we obtain
\[
\begin{aligned}
\sup_{t\in J}\int_{\R^d}
       (1+|\xi|^2)^s|\widehat\rho_t(\xi)|^2\dd\xi
 &\leq2^s|B_1|
       +5^sC_{n,J}\sum_{\ell=0}^\infty
                        2^{\ell(2s+d-4n/3)}\\
 &=2^s|B_1|
       +\frac{5^sC_{n,J}}{1-2^{2s+d-4n/3}}<\infty.
\end{aligned}
\]
The geometric series converges precisely because
$2s+d-4n/3<0$.
By Plancherel the inverse Fourier transform of
$\widehat\rho_t$ is an $H^s$ function. Its Fourier
transform agrees with that of the finite measure
$\rho_t$, so they are equal as tempered distributions.
Thus $\rho_t$ has an $H^s$ density, with the displayed
uniform bound. No density of the initial measure or
separate measure-to-function argument on the consensus
graph has been assumed.

We now prove continuity in $H^s$. Fix $t\in J$ and
any sequence $t_j\in J$ with $t_j\to t$. Choose $s'>s$
and an integer $n>\frac34(2s'+d)$. The preceding
calculation at order $s'$ gives
\(\sup_{\tau\in J}\int_{\R^d}(1+|\xi|^2)^{s'}
                |\widehat\rho_\tau(\xi)|^2\dd\xi\leq C\).

For each fixed $\xi$, the function
$v\mapsto e^{-i\xi\cdot v}$ is bounded and continuous.
Narrow continuity, applied to its real and imaginary
parts, therefore gives
$\widehat\rho_{t_j}(\xi)\to\widehat\rho_t(\xi)$.
On the fixed ball $|\xi|\leq L$, the integrand
$(1+|\xi|^2)^s|\widehat\rho_{t_j}-\widehat\rho_t|^2$
is bounded by $4(1+L^2)^s$, because each Fourier
transform has absolute value at most one. This bound
is integrable on that ball. Dominated convergence gives
\[
\lim_{j\to\infty}\int_{|\xi|\leq L}(1+|\xi|^2)^s
        |\widehat\rho_{t_j}(\xi)-\widehat\rho_t(\xi)|^2
                                                    \dd\xi=0.
\]
To control the complement uniformly in $j$, bound the
squared difference by twice the sum of the squared absolute
values and use
$(1+|\xi|^2)^{s-s'}\leq(1+L^2)^{s-s'}$ for $|\xi|>L$.
The uniform order-$s'$ bound yields
\[
\begin{aligned}
\int_{|\xi|>L}(1+|\xi|^2)^s
       |\widehat\rho_{t_j}-\widehat\rho_t|^2\dd\xi
 &\leq2(1+L^2)^{s-s'}\int_{|\xi|>L}(1+|\xi|^2)^{s'}
       \big(|\widehat\rho_{t_j}|^2+|\widehat\rho_t|^2\big)
                                                       \dd\xi\\
 &\leq4C(1+L^2)^{s-s'}.
\end{aligned}
\]
By Plancherel and the low-frequency limit, for every fixed
$L\geq1$ we consequently have
\(\limsup_{j\to\infty}\|\rho_{t_j}-\rho_t\|_{H^s}^2
 \leq\frac{4C}{(2\pi)^d}(1+L^2)^{s-s'}\).
First $j\to\infty$ was taken at fixed $L$; now let
$L\to\infty$. Since $s-s'<0$, the right side tends to
zero. Thus $\rho_{t_j}\to\rho_t$ in $H^s$. The time
$t$ and convergent sequence were arbitrary, proving
$\rho\in C(J;H^s)$ and hence \eqref{eq:all-sobolev-good}.
\end{proof}

\subsection{Propagation after the motion has ceased}
\label{sec:forward-propagation}

\begin{lemma}[Forward propagation of spatial Sobolev regularity]
\label[lemma]{lem:fourier-propagation}
Let $m\in C^1([0,T])$, and suppose the fourth-moment bound
\eqref{eq:fourth-moment} holds. If the given solution has
an $H^s$ density at time $t_*$, $s\geq0$, then for
$t_*\leq t\leq T$,
\begin{equation}\label{eq:fourier-propagation}
\|\rho_t\|_{H^s}^2
 \leq e^{C_s(t-t_*)}\|\rho_{t_*}\|_{H^s}^2.
\end{equation}
Here $C_s$ depends only on $d,\lambda,\sigma,s$; in particular it does not depend on $m$.
If $\rho_{t_*}\in H^s$ for every $s$, then
$\rho\in C([t_*,T];H^s)$ for every $s\geq0$.
No nonzero-velocity condition after $t_*$ is required.
\end{lemma}
This estimate is essentially contained in \cite{FS}, where it is derived in physical space rather than Fourier space, but only for integer $s$.
\begin{proof}
Fix an integer $N>s+d/2+4$, and use the integrable weight
\(W_\varepsilon(\xi) =\frac{(1+|\xi|^2)^s}{(1+\varepsilon|\xi|^2)^N},  0<\varepsilon\leq1\).
Its logarithm satisfies
\(|\xi|\,|\nabla_\xi\log W_\varepsilon| +|\xi|^2|\nabla_\xi^2\log W_\varepsilon|\leq C_{s,N,d}\),
as follows by differentiating the two logarithms
$s\log(1+|\xi|^2)-N\log(1+\varepsilon|\xi|^2)$.
Consequently, uniformly in $\varepsilon$,
\(\left|\Delta_\xi(|\xi|^2W_\varepsilon)\right| +\left|\divv_\xi(\xi W_\varepsilon)\right| \leq C_{s,N,d}W_\varepsilon\).
The moment bound makes $F$, $\nabla_\xi F$ and
$\nabla_\xi^2F$ bounded uniformly for $t\in[0,T]$
(with bounds allowed to depend on $T$), and
\eqref{eq:fourier-equation} gives
$|\partial_tF|\leq C(1+|\xi|^2)$. It is used
only to justify the integrations, not as the growth
constant in the energy estimate below.
For each fixed $\varepsilon$, the choice of $N$ therefore
justifies differentiation under the integral below.
Integration by parts, first with a large frequency cutoff
and then letting that cutoff tend to infinity, gives
\begin{equation}\label{eq:weighted-fourier-energy}
\begin{split}
\frac{\dd}{\dd t}\int W_\varepsilon|F|^2\dd\xi
={}&-\sigma^2\int W_\varepsilon|\xi|^2
                                      |\nabla_\xi F|^2\dd\xi\\
 &+\frac{\sigma^2}{2}\int
       \Delta_\xi(|\xi|^2W_\varepsilon)|F|^2\dd\xi
       +\lambda\int\divv_\xi(\xi W_\varepsilon)|F|^2\dd\xi\\
\leq{}&C\int W_\varepsilon|F|^2\dd\xi.
\end{split}
\end{equation}
Here $C$ depends only on $d,\lambda,\sigma,s,N$.
The imaginary potential again has zero real contribution.
All cutoff errors vanish by the polynomial decay of
$W_\varepsilon$ and the bounded derivatives just stated.
Gronwall's inequality, $W_\varepsilon\leq(1+|\xi|^2)^s$
at the initial time $t_*$, and monotone convergence as
$\varepsilon\downarrow0$ prove
\eqref{eq:fourier-propagation}. The harmless Fourier
normalization factor $(2\pi)^{-d}$ occurs on both sides.
Since $N$ is a function of $s$ and $d$, it is absorbed into $C_s$.

When all orders are available at $t_*$, apply the bound
also at an order $s'>s$, obtaining a uniform $H^{s'}$
bound on $[t_*,T]$. Fix $t_0\in[t_*,T]$ and split the
Fourier integral for $\|\rho_t-\rho_{t_0}\|_{H^s}^2$
at a radius $L\geq1$. On $|\xi|>L$, the inequality
$(1+|\xi|^2)^s\leq(1+L^2)^{s-s'}(1+|\xi|^2)^{s'}$
and the higher-order bound give
\[
\int_{|\xi|>L}(1+|\xi|^2)^s
 |\widehat\rho_t(\xi)-\widehat\rho_{t_0}(\xi)|^2\dd\xi
 \leq C(1+L^2)^{s-s'}.
\]
On $|\xi|\leq L$, narrow continuity gives pointwise
convergence of the Fourier transforms as $t\to t_0$;
the mass-one bound $|\widehat\rho_t|\leq1$ permits
dominated convergence on this bounded set. Thus the
low-frequency integral tends to zero for fixed $L$,
while the high-frequency bound tends to zero as
$L\to\infty$ because $s'<\infty$ and $s-s'<0$.
Taking these limits in that order and using Plancherel
proves continuity in $H^s$. The same argument applies
with one-sided limits at $t_*$ and $T$.
\end{proof}

\begin{proposition}[Any earlier nonconstant motion suffices]
\label[proposition]{prop:earlier-motion}
Let $m\in C^1([0,T])$ and $\rho_0\in\Pp_4$.
If $m$ is nonconstant on $[a,b]\subset[0,T)$, then
the solution belongs to $C([b,T];H^s)$ for every
$s\geq0$. If in addition $m\in C^k$,
then $\rho\in C^{\infty,k+1}(\R^d\times(b,T))$.
\end{proposition}

\begin{proof}
By the mean value theorem applied to a nonconstant
component, some $t_1\in(a,b)$ satisfies $m_j'(t_1)\neq0$.
Continuity supplies a smaller open interval
$I\Subset(a,b)$ and a constant $\delta>0$ satisfying
\eqref{eq:velocity-sign}. Choose any $t_*\in I$.
\Cref{prop:fourier-smoothing} gives all Sobolev orders
at $t_*$. \Cref{lem:fourier-propagation} carries them
to all later times, including a neighborhood of $b$.

Fix a spatial multi-index $\beta$ and choose
$s>|\beta|+d/2$. Sobolev embedding and the $H^s$
continuity just proved give
\[
\|\partial_v^\beta\rho_t-\partial_v^\beta\rho_{t_0}\|_\infty
 \leq C_{d,s,\beta}\|\rho_t-\rho_{t_0}\|_{H^s}
 \longrightarrow0\qquad(t\to t_0).
\]
For each fixed time, the spatial derivative is continuous
in $v$. Combining this spatial continuity with the uniform
convergence in the displayed estimate proves joint
continuity in $(v,t)$. Since $\beta$ was arbitrary,
this applies to every spatial derivative on
$\R^d\times(t_*,T)$.

Write the right-hand side of \eqref{eq:forward} as
$L_{m(t)}^*\rho$, with the same operator as in
\eqref{eq:time-bootstrap}. Its coefficients are continuous
in time and smooth in space. The preceding spatial
regularity therefore makes every
$\partial_v^\beta(L_{m(t)}^*\rho)$ jointly continuous.
Spatially differentiating the distributional equation
and integrating in time gives
\[
\partial_v^\beta\rho(v,t)-\partial_v^\beta\rho(v,t_0)
 =\int_{t_0}^t
       \partial_v^\beta\big(L_{m(\tau)}^*\rho(v,\tau)\big)
       \dd\tau.
\]
To justify the pointwise identity, first pair with a
compactly supported smooth spatial test. The resulting
continuous scalar function has the continuous
distributional derivative given by the right-hand side,
so the fundamental theorem of calculus applies. Both
sides are continuous in $v$, hence equality against all
such tests gives equality at every $v$.
It follows that $\partial_t\partial_v^\beta\rho$ exists
classically and equals the continuous integrand.

If $m\in C^k$, the coefficients of $L_{m(t)}^*$ have
continuous time derivatives through order $k$. Starting
from the identity just proved, induction and the product
rule give \eqref{eq:time-bootstrap} for $0\leq j\leq k$,
also after every spatial derivative. At step $j$, the
right-hand side uses only derivatives of $m$ through
order $j$ and already established time derivatives of
$\rho$ through order $j$, with arbitrary spatial
derivatives. It thus supplies the next continuous time
derivative. This proves $C^{\infty,k+1}$ on
$\R^d\times(t_*,T)$, and hence on the asserted strip.
The continuous densities represent every marginal by Fourier inversion.
\end{proof}

\subsection{The nonlinear bootstrap to a smooth consensus}
\label{sec:post-motion-bootstrap}

\begin{proposition}[The consensus becomes smooth after motion]
\label[proposition]{prop:post-motion-bootstrap}
Suppose $w>0$, $w,v_iw\in\mathcal X^2$, and
$\rho_0\in\Pp_4$. For the nonlinear CBO equation, if $m$
is nonconstant on $[a,b]\subset[0,T)$, then
\begin{equation}\label{eq:post-motion-infinity}
m\in C^\infty((b,T];\R^d),\qquad
\rho\in C^\infty(\R^d\times(b,T)).
\end{equation}
In particular no derivatives of $w$ beyond order two
are needed for this conclusion. Derivatives at $T$
are understood from the left.
\end{proposition}

\begin{proof}
\Cref{prop:observables}(ii) first gives $m\in C^1$.
By \cref{prop:earlier-motion}, on each interval
$J=[b',T]$ with $b<b'<T$ the density is uniformly bounded in
every spatial $H^s$, and hence all its spatial derivatives
are uniformly bounded on $\R^d\times J$.

Assume inductively that $m\in C^r((b,T])$, $r\geq1$,
with one-sided derivatives at $T$.
The proof of \cref{prop:earlier-motion} gives
$\rho\in C^{\infty,r+1}$ for $b<t<T$. This assertion
extends from the left to $T$: spatial Sobolev continuity
holds up to $T$, so the PDE right-hand side and the
successive right-hand sides in \eqref{eq:time-bootstrap}
extend continuously to $T$. Integrating each such
identity in time gives the corresponding one-sided
derivative at $T$, also after any spatial derivative.
We next obtain a spatial growth bound for these time
derivatives, which is needed before integrating them
against $w$ or $v_iw$ over all of $\R^d$.
The operator $L_{m(t)}^*$ in \eqref{eq:forward}
has spatial coefficients of degrees at most two, one,
and zero in its second-, first-, and zeroth-order terms.
For example, the first two time derivatives satisfy
$\partial_t\rho=L_{m(t)}^*\rho$ and
$\partial_t^2\rho=(L_{m(t)}^*)^2\rho+
(\partial_t L_{m(t)}^*)\rho$, where, for a fixed spatial test $\phi$,
\[
(\partial_t L_{m(t)}^*)\phi
 =-\sigma^2(v-m(t))\cdot m'(t)\,\Delta_v\phi
   -(\lambda+2\sigma^2)m'(t)\cdot\nabla_v\phi.
\]
Thus the second time derivative uses at most four spatial
derivatives, polynomial coefficients of degree at most
four, and only $m,m'$ as time coefficients.

More generally, for each $0\leq j\leq r+1$ we claim
\[
\partial_t^j\rho(v,t)
 =\sum_{|\beta|\leq2j}p_{j,\beta}(v,t)
                              \partial_v^\beta\rho(v,t),
\]
where $p_{j,\beta}$ is a polynomial in $v$ of degree at
most $2j$, whose coefficients are polynomials in the
components of $m,\ldots,m^{(j-1)}$ when $j\geq1$.
The assertion at $j=0$ holds with $p_{0,0}=1$.
Suppose it holds at some $j\leq r$. Differentiating gives
\[
\partial_t^{j+1}\rho
 =\sum_{|\beta|\leq2j}
  \left((\partial_t p_{j,\beta})\partial_v^\beta\rho
        +p_{j,\beta}\partial_v^\beta(L_{m(t)}^*\rho)\right).
\]
The first term does not increase the polynomial degree
in $v$ and uses at most $m^{(j)}$.
For the second term, the spatial Leibniz rule gives
\[
\begin{aligned}
\partial_v^\beta(L_{m(t)}^*\rho)
={}&\frac{\sigma^2}{2}\sum_{\gamma\leq\beta}\binom\beta\gamma
  \big(\partial_v^\gamma|v-m(t)|^2\big)
                    \Delta_v\partial_v^{\beta-\gamma}\rho\\
 &+(\lambda+2\sigma^2)
  \sum_{\gamma\leq\beta}\binom\beta\gamma
   \big(\partial_v^\gamma(v-m(t))\big)\cdot
                    \nabla_v\partial_v^{\beta-\gamma}\rho
   +d(\lambda+\sigma^2)\partial_v^\beta\rho.
\end{aligned}
\]
Terms in the two sums vanish for $|\gamma|>2$ and
$|\gamma|>1$, respectively. Multiplication by $p_{j,\beta}$ adds at
most two polynomial degrees, and the density is
differentiated at most $|\beta|+2\leq2j+2$ times.
This proves the asserted representation at $j+1$,
including its dependence on derivatives of $m$.

Since those derivatives of $m$ are bounded on the fixed
compact interval $J$, the representation gives
\[
\begin{aligned}
\sup_{t\in J}|\partial_t^j\rho(v,t)|
 &\leq C(1+|v|)^{2j}
       \sum_{|\beta|\leq2j}
                    \sup_{t\in J}\|\partial_v^\beta\rho_t\|_\infty\\
 &\leq C(1+|v|)^{2j}
       \sup_{t\in J}\|\rho_t\|_{H^{2j+\lfloor d/2\rfloor+1}}.
\end{aligned}
\]
The last inequality is Sobolev embedding, since
$2j+\lfloor d/2\rfloor+1>2j+d/2$.
Here $C$ depends only on $d,j,\lambda,\sigma$ and, for $j\geq1$, the supremum norms of $m,\ldots,m^{(j-1)}$ on $J$. Absorbing that finite norm
into the constant proves
\begin{equation}\label{eq:polynomial-time-bounds}
\sup_{t\in J}|\partial_t^j\rho(v,t)|
       \leq C_{j,J}(1+|v|)^{2j},
\qquad 0\leq j\leq r+1.
\end{equation}
Here $C_{j,J}$ depends only on $d,j,\lambda,\sigma$,
the uniform $H^{2j+\lfloor d/2\rfloor+1}$ norm of
$\rho$ on $J$, and, for $j\geq1$, the supremum norms
of $m,\ldots,m^{(j-1)}$ on $J$.
For each weight $\psi=w$ or $\psi=v_iw$, rapid decay gives
$\int |\psi(v)|(1+|v|)^{2(r+1)}\dd v<\infty$.
Consequently \eqref{eq:polynomial-time-bounds} supplies
an integrable majorant, independent of $t\in J$, for
every $\psi\partial_t^j\rho$, $0\leq j\leq r+1$.
For $0\leq j\leq r$, multiply the pointwise fundamental
theorem of calculus for $\partial_t^j\rho$ by $\psi$ and
use this majorant to apply Fubini. For $s,t\in J$ we obtain
\[
\int \psi(v)\partial_t^j\rho(v,t)\dd v
 -\int \psi(v)\partial_t^j\rho(v,s)\dd v
 =\int_s^t\!\int \psi(v)\partial_\tau^{j+1}\rho(v,\tau)
                                               \dd v\dd\tau.
\]
The inner integral is continuous in $\tau$ by dominated
convergence and the pointwise continuity already proved.
Thus differentiation of this scalar identity, successively
in $j$, gives
\(\frac{\dd^j}{\dd t^j}\int \psi(v)\rho(v,t)\dd v =\int \psi(v)\partial_t^j\rho(v,t)\dd v,  0\leq j\leq r+1\).
No derivative falls on the observable weight $\psi$;
in particular no higher derivatives of $f$ are required.
Hence $A(t)=\int vw\dd\rho_t$ and
$B(t)=\int w\dd\rho_t$ belong to $C^{r+1}$ on $(b,T]$,
including the one-sided derivatives at $T$.
On $J$, continuity and positivity give $\min_J B>0$,
so the quotient $m=A/B$ is $C^{r+1}$ there, including
its one-sided derivatives at $T$. Since $b'>b$ was
arbitrary, $m\in C^{r+1}((b,T])$.
This closes the induction starting from $r=1$.
The time bootstrap for the PDE then gives all time
derivatives of the density as well, proving
\eqref{eq:post-motion-infinity}.
\end{proof}

\subsection{Completion of the proof}
\label{sec:onset-completion}

\begin{proof}[Proof of \cref{thm:main}]
Let $\rho=(\rho_t)_{t\geq0}$ be the unique global weak solution
given by \cref{thm:global-existence}.
Fix $T>0$. By \cref{prop:observables}(ii),
$m\in C^k([0,T];\R^d)$, in particular $m\in C^1([0,T];\R^d)$.
Let $t_c$ be defined by \cref{def:onset} on this interval.

If $t_c>0$, then $m(t)=m(0)$ on $[0,t_c]$.
Apply \cref{lem:off-center} on $(0,t_c)$ with arbitrarily large $r$.
The constant trajectory satisfies every such $C^r$ assumption,
so $\rho$ is jointly smooth on
$(\R^d\setminus\{m(0)\})\times(0,t_c)$. This proves \textup{(i)}.

If $t_c<T$, fix $a\in(t_c,T)$ and choose $\tau\in(t_c,a)$.
By \cref{def:onset}, $m$ is nonconstant on $[0,\tau]$.
\Cref{prop:earlier-motion} therefore gives
$\rho\in C([a,T];H^s(\R^d))$ for every $s\geq0$.
\Cref{prop:post-motion-bootstrap} gives
$m\in C^\infty((\tau,T];\R^d)$ and
$\rho\in C^\infty(\R^d\times(\tau,T))$.
Since $a>t_c$ is arbitrary, this proves \textup{(ii)}.
Together with the constant trajectory before $t_c$, it also proves
\eqref{eq:consensus-onset-regularity}.

If $0<t_c<T$, all derivatives of $m$ through order $k$ vanish
on $(0,t_c)$. Their continuity at $t_c$ gives
$m^{(j)}(t_c)=0$ for $1\leq j\leq k$.
For $v_0\neq m(t_c)$, continuity of $m$ provides a cylinder around
$(v_0,t_c)$ that stays away from the consensus graph.
Applying \cref{lem:off-center} with $r=k$ on this cylinder gives
$\rho\in C^{\infty,k+1}$ across the transition, proving \textup{(iii)}.

When $t_c=0$, the argument for \textup{(ii)} applies to every positive time.
When $t_c=T$, $m$ is constant on $[0,T]$ and hence smooth there,
and \textup{(i)} applies on $(0,T)$.
Since $T>0$ was arbitrary, the conclusions hold on every finite
interval of the global solution.
\end{proof}

\section{Positivity and full support}
\label{sec:positivity}

Fornasier--Sun \cite[Section~4.2 and Theorem~4.5, pp.~28--30]{FS}
prove positivity away from the consensus point in dimension $d>1$
for classical solutions starting from initial densities in every $H^m$.
The following theorem extends positivity and full support to non-Dirac
measure initial data in every dimension, including dimension one.

\begin{theorem}[Positivity and full support from rough initial data]
\label[theorem]{thm:fs-positivity}
Assume \cref{ass:main} with $k=1$, let $d\geq1$, and
let $\rho$ be a weak solution with initial law $\rho_0$
which is not a Dirac measure. Then, we have 
\begin{equation}\label{eq:fs-positivity}
\rho(v,t)>0\qquad(0<t<T,\ v\neq m(t)).
\end{equation}
Consequently,
\begin{equation}\label{eq:fs-full-support}
\supp\rho_t=\R^d,\qquad
\rho_t(B_r(v))>0
\quad(0<t<T,\ v\in\R^d,\ r>0).
\end{equation}
\end{theorem}

\begin{proof}
By \cref{prop:observables}(ii), $m\in C^1([0,T])$.
We first exclude concentration at a single point at any finite time.
Write $\bar v_t=\int v\dd\rho_t(v)$ and
$V_t=\int|v-\bar v_t|^2\dd\rho_t(v)$.
Testing the weak equation with $v$ and $|v|^2$, using spatial
cutoffs justified by the fourth-moment bound in
\cref{lem:fourth-moment}, gives
$\bar v_t'=-\lambda(\bar v_t-m(t))$ and
$V_t'=(-2\lambda+d\sigma^2)V_t+d\sigma^2|\bar v_t-m(t)|^2$.
Since $\rho_0$ is non-Dirac, $V_0>0$, and hence
$V_t\geq e^{(-2\lambda+d\sigma^2)t}V_0>0$.
Thus $\rho_t$ is non-Dirac for every $t\in[0,T]$.

In moving coordinates,
put $y=v-m(t)$ and $\widetilde\rho(y,t)=\rho(y+m(t),t)$ for $y\neq0$.
The local regularity from \cref{lem:off-center} justifies
the classical chain rule, which gives
\begin{equation}\label{eq:fs-moving-positive}
\partial_t \widetilde\rho
 =\frac{\sigma^2}{2}|y|^2\Delta_y \widetilde\rho
       +\big((\lambda+2\sigma^2)y+m'(t)\big)\cdot\nabla_y \widetilde\rho
       +d(\lambda+\sigma^2)\widetilde\rho.
\end{equation}
Equivalently,
\(\partial_t \widetilde\rho =\divv_y\left(\frac{\sigma^2}{2}|y|^2\nabla_y \widetilde\rho +\big((\lambda+\sigma^2)y+m'(t)\big)\widetilde\rho\right)\).
Indeed, $\nabla(\sigma^2|y|^2/2)=\sigma^2y$ and
$\divv_y((\lambda+\sigma^2)y+m'(t))=d(\lambda+\sigma^2)$,
so expansion
recovers \eqref{eq:fs-moving-positive}.
The term $m'(t)\cdot\nabla \widetilde\rho$ is exactly
$\divv(m'(t)\widetilde\rho)$; no time derivative of $m'$ occurs.

On any annulus $A=\{y:r_-<|y|<r_+\}$ with
$0<r_-<r_+<\infty$ and any open interval $J\Subset(0,T)$,
the following bounds hold for $y\in A$, $t\in J$ and $z\in\R^d$:
\(\frac{\sigma^2r_-^2}{2}|z|^2 \leq\frac{\sigma^2}{2}|y|^2|z|^2 \leq\frac{\sigma^2r_+^2}{2}|z|^2\), and \(|(\lambda+\sigma^2)y+m'(t)|\leq(\lambda+\sigma^2)r_+ +\|m'\|_{L^\infty(J)}\).
These bounds are precisely sufficient for the divergence
form parabolic Harnack inequality in
\cite[Theorem~H, p.~618]{Aronson1968}.
In the notation of that paper, take
\(A_{ij}=\frac{\sigma^2}{2}|y|^2\delta_{ij}, A_j=(\lambda+\sigma^2)y_j+m_j'(t),  B_j=C=0\).
The coefficients are measurable and bounded, and the
principal matrix is uniformly positive definite.
Thus hypotheses \textup{(H.1)--(H.3)} there hold on each
relatively compact cylinder; see
\cite[pp.~608--609]{Aronson1968}, including the stated
local extension to globally bounded coefficients.
In particular the theorem does not assume smoothness
in time of the drift coefficient. By \cref{lem:off-center}, $\widetilde\rho$ is a
nonnegative classical solution on these cylinders, hence
also a local energy weak solution, as required by that
theorem. Its continuous representative permits pointwise
application of the inequality.

Fix a time $t\in(0,T)$ and a spatial point $v\in\R^d$ with
$v\neq m(t)$. We prove that $\rho(v,t)>0$.
Choose any earlier time $s\in(0,t)$.
If $d\geq2$, then
$\rho_s(\R^d\setminus\{m(s)\})>0$, so
\cref{lem:off-center} supplies a point $v_1\neq m(s)$ with
$\rho(v_1,s)>0$.
If $d=1$, the identity
$\int_{\R}(u-m(s))w(u)\dd\rho_s(u)=0$ and $w>0$ imply
$\rho_s((-\infty,m(s)))>0$ and
$\rho_s((m(s),\infty))>0$.
Indeed, if either half-line had zero mass, the integrand would
have one sign, and its zero integral would force
$\rho_s=\delta_{m(s)}$, a contradiction.
The density is continuous on each half-line, so each half-line contains
a point where the density is positive.
For any spatial point $v\in\R$ with $v\neq m(t)$, choose $v_1$ as follows:
If $v>m(t)$, choose $v_1>m(s)$ with $\rho(v_1,s)>0$.
If $v<m(t)$, choose $v_1<m(s)$ with $\rho(v_1,s)>0$.
In both cases, $(v_1-m(s))(v-m(t))>0$.

Put $y_1=v_1-m(s)$ and $y_2=v-m(t)$.
Choose $0<r_-<r_+$ so that $A=\{y:r_-<|y|<r_+\}$ contains
both points in its interior, and choose an open interval
$J\Subset(0,T)$ containing $[s,t]$.
The points $y_1,y_2$ lie in the same connected component of $A$:
for $d\geq2$ the annulus is connected, while for $d=1$ they
have the same sign because $y_1y_2=(v_1-m(s))(v-m(t))>0$.
Choose a continuous path
$\gamma:[0,1]\to A$ joining $y_1$ to $y_2$.
Its compact image has positive distance from $\partial A$.
Choose $\ell>0$ such that every ball of radius $2\ell$
centered on that image is contained in $A$. By uniform
continuity choose an integer $N$ so large that
$|\gamma((j+1)/N)-\gamma(j/N)|<\ell/2$ for every $j$.
Set $z_j=\gamma(j/N)$ and
$t_j=s+j(t-s)/N$, $j=0,\ldots,N$.
For each adjacent pair apply
\cite[Theorem~H, p.~618]{Aronson1968} on
$B_{2\ell}(z_j)\times J$ with the inner convex ball
$B_\ell(z_j)$, which contains both $z_j$ and $z_{j+1}$.
It gives a finite constant $C_j$ such that
\(\widetilde\rho(z_j,t_j)\leq C_j \widetilde\rho(z_{j+1},t_{j+1})\).
All times are interior to $J$, and $t_{j+1}>t_j$;
each $C_j$ depends only on $d,\lambda,\sigma,r_-,r_+$, $\|m'\|_{L^\infty(J)}$, the chosen balls, $J$ and the two comparison points, not on the solution $\widetilde\rho$.
A finite iteration gives
$0<\widetilde\rho(y_1,s)\leq(\prod_{j=0}^{N-1}C_j)\widetilde\rho(y_2,t)$.
Hence $\rho(v,t)>0$, proving
\eqref{eq:fs-positivity}.

Every nonempty open ball contains a smaller ball whose
closure avoids $m(t)$. The positive continuous density
on that smaller ball has positive integral. Thus every
open ball has positive $\rho_t$ mass, proving
\eqref{eq:fs-full-support}, including balls centered at
$m(t)$. A possible singular component at the center
does not affect this support conclusion.
\end{proof}

The full-support result gives an application to global optimization,
following Fornasier--Sun \cite[Introduction, equations~(7)--(9)]{FS}.
Assume additionally that $f$ satisfies \textup{(A1)--(A2)} of
\cite[Definition~3.5]{FKR2024}, with global minimizer $v^*$, and that
$2\lambda>d\sigma^2$. Starting from any non-Dirac
$\rho_0\in\Pp_4$, run \eqref{eq:cbo} with a fixed parameter until
any time $\tau>0$. The law $\mu=\rho_\tau$ has full support by
\cref{thm:fs-positivity}, so $\mu(B_r(v^*))>0$ for every $r>0$.
For any prescribed accuracy $0<\varepsilon<W_2^2(\mu,\delta_{v^*})$
and $\theta\in(0,1)$,
choose a sufficiently large optimization parameter $\alpha$ after
fixing $\mu$, and restart from $\widetilde\rho_0=\mu$.
Then \cite[Theorem~3.7]{FKR2024} gives
\begin{equation}\label{eq:full-support-optimization}
W_2^2(\widetilde\rho_s,\delta_{v^*})
\leq W_2^2(\mu,\delta_{v^*})
       e^{-(1-\theta)(2\lambda-d\sigma^2)s}
\end{equation}
until the squared distance reaches $\varepsilon$ at some finite time
$T$. Thus \eqref{eq:full-support-optimization} yields global
optimization to any prescribed accuracy without requiring the original
initial measure to have mass near the unknown minimizer.

This also leads to convergence of the discrete CBO algorithm.
Under the additional objective assumptions of
\cite[Proposition~3.11]{FKR2024}, initialize $N$ particles
$V_0^1,\ldots,V_0^N$ independently with law $\mu$ and, for
$i=1,\ldots,N$ and $k=0,\ldots,K-1$, use
\begin{equation}\label{eq:discrete-cbo}
V_{k+1}^i=V_k^i-\lambda h(V_k^i-m_k^N)
       +\sigma\sqrt h\,|V_k^i-m_k^N|Z_k^i.
\end{equation}
Here $h=T/K$, $T$ is the terminal time above,
$m_k^N=\sum_{j=1}^N V_k^j e^{-\alpha f(V_k^j)}
/\sum_{j=1}^N e^{-\alpha f(V_k^j)}$, and the $Z_k^i$ are independent
standard Gaussian vectors in $\R^d$, independent of the initial particles.
For any total squared-error tolerance $\varepsilon_{\mathrm{alg}}>0$
and failure probability $\delta\in(0,1)$, choose the continuous error
$\varepsilon$ sufficiently small, then $\alpha$ and $N$ sufficiently
large, and finally $h$ sufficiently small.
Then \cite[Theorem~3.8, with Euler--Maruyama order $1/2$]{FKR2024}
gives
$\mathbb P\big(\big|N^{-1}\sum_{i=1}^N V_K^i-v^*\big|^2
\leq\varepsilon_{\mathrm{alg}}\big)\geq1-\delta$
for the algorithm \eqref{eq:discrete-cbo}.

\section*{Statement on the use of AI}
The author used AI tools to assist with literature searches, improve the wording of the manuscript, and simplify some calculations. The original proof in Section~3 worked directly with $\widehat{\rho}_t$, whose evolution equation contains additional terms and leads to more involved estimates. When asked to check this proof, an AI tool suggested estimating the centered Fourier transform $F$ instead, which simplified the argument. AI tools were also used to carry out the routine but lengthy calculations for the explicit example in Appendix~B. The author has reviewed and checked the entire manuscript and takes full responsibility for its contents.

\appendix

\section{Global well-posedness of the weak solution in Theorem~\ref{thm:main}}
\label{sec:existence}

This appendix proves global well-posedness under the weight
assumptions of \cref{thm:main}. We first construct local solutions by Gaussian transition kernels and establish a second-moment
continuation criterion. We then exclude a finite maximal endpoint:
a short-time lower bound on the linear transition law bounds the
consensus denominator from below, and rapid decay of the weight closes
a Gronwall inequality for the accumulated forcing.
The latter argument uses the Gaussian perturbation estimate of
\cite[Proposition~3.1(i)]{BallyCaramellino2013}, with explicit bounds
on the first two Malliavin derivatives of a linear SDE.
It requires neither an initial density nor
the regularization theorem. Finally, positivity of the denominator
and a synchronous coupling give $W_2$ stability and uniqueness in
\cref{prop:weak-uniqueness}. This last argument uses no Malliavin
derivatives. The regularity and positivity proofs in the main text
remain entirely analytic.

\subsection{Local existence and continuation}

Throughout this appendix assume
\begin{equation}\label{eq:existence-weight}
w>0,\qquad w,v_1w,\ldots,v_dw\in\mathcal X^2,
\qquad \rho_0\in\Pp_4(\R^d).
\end{equation}
These hypotheses follow from \cref{ass:main} for every $k\geq1$.
Write, for a probability measure $\mu$,
\(B(\mu)=\int w\dd\mu, A(\mu)=\int vw\dd\mu\).
The consensus quotient $A(\mu)/B(\mu)$ is defined for every probability measure:
the numerator is bounded and $B(\mu)>0$. Fix
$\varepsilon>0$ and define
\(\mathcal M_{w,\varepsilon}(\mu) :=\frac{A(\mu)}{B(\mu)\vee\varepsilon}\).
Then $|\mathcal M_{w,\varepsilon}(\mu)|\leq\|vw\|_\infty/\varepsilon$.
Both $w$ and the vector-valued function $v\mapsto vw(v)$
are bounded and globally Lipschitz. A coupling of $\mu,\nu\in\Pp_1$
gives
\begin{equation}\label{eq:cutoff-consensus-lipschitz}
|\mathcal M_{w,\varepsilon}(\mu)
       -\mathcal M_{w,\varepsilon}(\nu)|
\leq \left(\frac{\|\nabla(vw)\|_\infty}{\varepsilon}
       +\frac{\|vw\|_\infty\|\nabla w\|_\infty}{\varepsilon^2}\right)
       W_1(\mu,\nu).
\end{equation}
Indeed, $|A(\mu)-A(\nu)|\leq\|\nabla(vw)\|_\infty W_1(\mu,\nu)$ and
$|B(\mu)-B(\nu)|\leq\|\nabla w\|_\infty W_1(\mu,\nu)$, while
$b\mapsto b\vee\varepsilon$ is $1$-Lipschitz and is
bounded below by $\varepsilon$.

\begin{lemma}[Existence for the bounded consensus approximation]
\label[lemma]{lem:cutoff-existence}
For every $\varepsilon>0$ and $T>0$, there exists
$\rho^\varepsilon\in C([0,T];\Pp_2(\R^d))$, continuous
in $W_2$, with initial value $\rho_0$, satisfying
\begin{equation}\label{eq:cutoff-fpk}
\partial_t\rho^\varepsilon
 =\lambda\divv\big((v-m_\varepsilon(t))\rho^\varepsilon\big)
       +\frac{\sigma^2}{2}\Delta
          \big(|v-m_\varepsilon(t)|^2\rho^\varepsilon\big),
\quad
m_\varepsilon(t)=\mathcal M_{w,\varepsilon}(\rho_t^\varepsilon),
\end{equation}
in the weak sense. Moreover,
\begin{equation}\label{eq:cutoff-fourth-bound}
\sup_{0\leq t\leq T}\int|v|^4\dd\rho_t^\varepsilon(v)
 \leq C_{\varepsilon,T}\left(1+\int|v|^4\dd\rho_0(v)\right).
\end{equation}
Here $C_{\varepsilon,T}$ depends only on $d,\lambda,\sigma,\varepsilon,T$ and $\|vw\|_\infty$. Uniformity as
$\varepsilon\downarrow0$ is not asserted.
\end{lemma}

\begin{proof}
In this proof $C$ depends only on $d,\lambda,\sigma$.
The constants $C_{\varepsilon,T}$ may additionally depend on $\varepsilon,T,\|vw\|_\infty$ and $\int|v|^4\dd\rho_0$, but not on the mesh size $h$, the number of steps $N$ or the step index $n$.
Here the initial moment is
included in $C_{\varepsilon,T}$; in the statement it was
factored out explicitly. The test-dependent constant
$C_\phi$ and remainder constants are specified in Step~5.

\emph{Step 1: positive discrete approximations.}
Let $h=T/N\leq\min\{1,(2\lambda)^{-1}\}$ and
$t_n=nh$. Denote by $\gamma_d$ the standard Gaussian
probability measure on $\R^d$. Starting from
$\mu_0^h=\rho_0$, put
$q_n=\mathcal M_{w,\varepsilon}(\mu_n^h)$ and define
\begin{equation}\label{eq:gaussian-euler-scheme}
\begin{split}
\int\phi\dd\mu_{n+1}^h
 =\int\!\int\phi\big(&v-\lambda h(v-q_n)
             +\sigma\sqrt h\,|v-q_n|z\big)
             \dd\gamma_d(z)\dd\mu_n^h(v).
\end{split}
\end{equation}
Here \eqref{eq:gaussian-euler-scheme} defines $\mu_{n+1}^h$
by requiring the identity to hold for every bounded continuous
function $\phi$. Explicitly, take independent random variables
$X_n\sim\mu_n^h$ and $Z_{n+1}\sim N(0,I_d)$; then $\mu_{n+1}^h$
is the law of
\(X_n-\lambda h(X_n-q_n)+\sigma\sqrt h\,|X_n-q_n|Z_{n+1}\).
The right-hand side of \eqref{eq:gaussian-euler-scheme} is precisely
the expectation of $\phi$ evaluated at this random variable.
Thus $\mu_{n+1}^h$ is a uniquely defined probability measure,
even when $v=q_n$, where the transition is a Dirac measure.
In particular positivity and total mass one are exact.

\emph{Step 2: a uniform fourth-moment bound.}
For $a\in\R^d$ and $b\in\R$, we have \(\int|a+bz|^4\dd\gamma_d
=|a|^4+2(d+2)b^2|a|^2+d(d+2)b^4\).
Here $a=(1-\lambda h)v+\lambda hq_n$ and
$b=\sigma\sqrt h\,|v-q_n|$.
Since $0\leq\lambda h\leq1$, convexity gives
$|a|^4\leq(1-\lambda h)|v|^4+\lambda h|q_n|^4$
and $|a|^2\leq(1-\lambda h)|v|^2+\lambda h|q_n|^2
\leq |v|^2+|q_n|^2$. Consequently,
$|v-q_n|^2|a|^2\leq2(|v|^2+|q_n|^2)^2
\leq4(|v|^4+|q_n|^4)$ and
$|v-q_n|^4\leq4(|v|^2+|q_n|^2)^2
\leq8(|v|^4+|q_n|^4)$.
Substitution in the Gaussian identity, using $h^2\leq h$, yields
\[
\begin{aligned}
\int|a+bz|^4\dd\gamma_d(z)
&\leq(1-\lambda h)|v|^4+\lambda h|q_n|^4\\
&\quad+8(d+2)\sigma^2h(|v|^4+|q_n|^4)
       +8d(d+2)\sigma^4h^2(|v|^4+|q_n|^4)\\
&\leq(1+Ch)|v|^4+Ch|q_n|^4.
\end{aligned}
\]
Here one may take $C=\lambda+8(d+2)\sigma^2+8d(d+2)\sigma^4$.
Integrating against $\mu_n^h$ in \eqref{eq:gaussian-euler-scheme}
therefore gives the fourth-moment recursion
\(\int|v|^4\dd\mu_{n+1}^h \leq(1+Ch)\int|v|^4\dd\mu_n^h +Ch\left(\frac{\|vw\|_\infty}{\varepsilon}\right)^4.\)
Starting from $\mu_0^h=\rho_0$ and iterating this inequality,
\(\int|v|^4\dd\mu_n^h\leq e^{CT}\left[\int|v|^4\dd\rho_0 +CT\left(\frac{\|vw\|_\infty}{\varepsilon}\right)^4\right], 0\leq n\leq N,\)
since $nh\leq T$.
Because $\mu_n^h$ is a probability measure, Cauchy--Schwarz gives
$\int|v|^2\dd\mu_n^h\leq(\int|v|^4\dd\mu_n^h)^{1/2}$.
Thus the fourth-moment bound also supplies the second-moment bound:
\begin{equation}\label{eq:uniform-discrete-moments}
\max_{0\leq n\leq N}\int(|v|^2+|v|^4)\dd\mu_n^h
 \leq C_{\varepsilon,T}.
\end{equation}

\emph{Step 3: equicontinuity in $W_2$.}
Iterating the kernels in \eqref{eq:gaussian-euler-scheme}
defines a joint measure on the finite path
$(X_0,\ldots,X_N)$, with marginals $\mu_n^h$.
Equivalently, its integrals are successive integrations
against $\rho_0$ and $\gamma_d^{\otimes N}$, with
\(X_{n+1}-X_n=-\lambda h(X_n-q_n) +\sigma\sqrt h\,|X_n-q_n|Z_{n+1}\).
The Gaussian variable $Z_{n+1}$ is independent of the
preceding coordinates, has mean zero and second moment
$d$. Thus the increments $\sigma\sqrt h\,|X_n-q_n|Z_{n+1}$
at different steps are orthogonal in
$L^2$. With expectation denoting integration against
this finite joint measure, for $j<l$ we obtain
\begin{align*}
\mathbb E|X_l-X_j|^2
&\leq 2\lambda^2(t_l-t_j)
        \sum_{n=j}^{l-1}h\mathbb E|X_n-q_n|^2
       +2d\sigma^2\sum_{n=j}^{l-1}h\mathbb E|X_n-q_n|^2\\
&\leq C_{\varepsilon,T}
       \big((t_l-t_j)^2+(t_l-t_j)\big)
 \leq C_{\varepsilon,T}(t_l-t_j).
\end{align*}
This coupling proves the corresponding $W_2^2$ bound
at grid times. For $t=t_n+\theta h$, $0\leq\theta\leq1$,
interpolate by
\(\mu_t^h=(1-\theta)\mu_n^h+\theta\mu_{n+1}^h\).
For $s,t$ in the same grid cell, a coupling that leaves
the common mixture mass fixed and transports only the
remaining mass gives
\(W_2^2(\mu_t^h,\mu_s^h) \leq\frac{|t-s|}{h}\mathbb E|X_{n+1}-X_n|^2 \leq C_{\varepsilon,T}|t-s|\).
If $|t-s|<h$ and the points lie in adjacent cells, the
triangle inequality gives the same estimate with an
additional fixed factor. If $|t-s|\geq h$, compare with
the two neighboring grid times and use the grid estimate.
In every case,
\begin{equation}\label{eq:discrete-W2-modulus}
W_2(\mu_t^h,\mu_s^h)
 \leq C_{\varepsilon,T}|t-s|^{1/2},\qquad s,t\in[0,T].
\end{equation}
No uniformity as $\varepsilon\downarrow0$ is asserted.

\emph{Step 4: compactness of the curves.}
By Steps~2--3, the interpolated curves $\mu^h$ are uniformly
bounded and equicontinuous in $C([0,T];\Pp_2(\R^d))$ with respect
to $W_2$. More precisely, their values belong to the compact metric space
\(K_{\varepsilon,T} :=\left\{\mu\in\Pp_2(\R^d): \int|v|^4\dd\mu\leq C_{\varepsilon,T}\right\}\),
equipped with the distance $W_2$,
where the fourth-moment bound extends to the interpolation by convexity.
Apply the metric Arzel\`a--Ascoli theorem, using the common modulus
in \eqref{eq:discrete-W2-modulus}, in the space
$C([0,T];K_{\varepsilon,T})$ equipped with the uniform distance
$\sup_{t\in[0,T]}W_2(\mu_t,\nu_t)$.
There are a subsequence $h_j\to0$ and
$\rho^\varepsilon\in C([0,T];K_{\varepsilon,T})$ such that
\begin{equation}\label{eq:cutoff-curve-convergence}
\sup_{t\in[0,T]}W_2(\mu_t^{h_j},\rho_t^\varepsilon)
 \longrightarrow0.
\end{equation}
Since $\mu_0^h=\rho_0$, the limit has initial value $\rho_0$.
Lower semicontinuity of the fourth moment and the bound in Step~2
give \eqref{eq:cutoff-fourth-bound}. We suppress the subsequence index below.

\emph{Step 5: consistency and the weak equation.}
For $\phi\in C_c^\infty(\R^d)$, expand $\phi(v+\Delta)$
to second order, where
$\Delta=-\lambda h(v-q_n)+\sigma\sqrt h\,|v-q_n|z$.
Gaussian integration gives
\[
\int\Delta\dd\gamma_d=-\lambda h(v-q_n),\qquad
\int\Delta\Delta^{\mathsf T}\dd\gamma_d
 =\sigma^2h|v-q_n|^2I_d
       +\lambda^2h^2(v-q_n)(v-q_n)^{\mathsf T}.
\]
The Taylor remainder is bounded by $C_\phi|\Delta|^3$,
where $C_\phi$ depends only on $d,\|\nabla^3\phi\|_\infty$.
Since $\int|\Delta|^3\dd\gamma_d
\leq Ch^{3/2}(|v|^3+|q_n|^3)$, Step~2 yields
\(\int\phi\dd\mu_{n+1}^h-\int\phi\dd\mu_n^h =h\int L_{q_n}\phi\dd\mu_n^h +O_{\phi,\varepsilon,T}(h^{3/2}).\)
The $O(h^2)$ drift contribution is included in this remainder.
The implicit constant depends on $\phi$ only through $\|\phi\|_{C^3}$.
Summing over the time steps gives a total error $O(\sqrt h)$.

Let $\bar\mu_t^h=\mu_n^h$ for $t\in[t_n,t_{n+1})$.
By \eqref{eq:discrete-W2-modulus}--\eqref{eq:cutoff-curve-convergence},
$\bar\mu^h\to\rho^\varepsilon$ uniformly in $W_2$;
\eqref{eq:cutoff-consensus-lipschitz} then gives uniform convergence
of the stepwise $q_n$ to
$m_\varepsilon(t)=\mathcal M_{w,\varepsilon}(\rho_t^\varepsilon)$.
Since $L_q\phi$ has compact support and depends continuously on
bounded $q$ in the uniform norm, the summed identity passes to
\(\int\phi\dd\rho_t^\varepsilon-\int\phi\dd\rho_0 =\int_0^t\int L_{m_\varepsilon(s)}\phi\dd\rho_s^\varepsilon\dd s, t\in[0,T].\)
Integration against smooth time tests, followed by approximation
by finite sums of product tests, gives \eqref{eq:cutoff-fpk}
for all compactly supported smooth space--time tests.
\end{proof}

\begin{theorem}[Local existence under the original assumptions]
\label[theorem]{thm:local-existence}
Under \eqref{eq:existence-weight}, there is a $T_0>0$
and a solution of \eqref{eq:cbo} on $[0,T_0]$ satisfying
\(\rho\in C([0,T_0];\Pp_2(\R^d)), \rho_{t=0}=\rho_0, \sup_{0\leq t\leq T_0}\int|v|^4\dd\rho_t<\infty\).
In particular the solution is narrowly continuous.
If the higher-order weight condition in \cref{ass:main}
holds, its actual consensus satisfies $m\in C^k([0,T_0])$.
\end{theorem}

\begin{proof}
Choose $\varepsilon=B(\rho_0)/4>0$.
Apply \cref{lem:cutoff-existence}
on $[0,1]$. The bounded, continuous trajectory
$m_\varepsilon$ has $|m_\varepsilon|\leq\|vw\|_\infty/\varepsilon$.
Since $w\in\mathcal X^2$, \cref{lem:observable} applies
to this already constructed solution of the linear equation
with prescribed trajectory $m_\varepsilon$. It gives
\(\frac{\dd}{\dd t}\int w\dd\rho_t^\varepsilon =\int L_{m_\varepsilon(t)}w\dd\rho_t^\varepsilon\).
All norms in the following number are finite:
\(D_\varepsilon={}\lambda\|v\cdot\nabla w\|_\infty +\frac{\lambda\|vw\|_\infty}{\varepsilon}\|\nabla w\|_\infty +\sigma^2\||v|^2\Delta w\|_\infty +\frac{\sigma^2\|vw\|_\infty^2}{\varepsilon^2}\|\Delta w\|_\infty\).
The generator formula and $|v-m|^2\leq2|v|^2+2|m|^2$
show that $\|L_{m_\varepsilon(t)}w\|_\infty\leq D_\varepsilon$.
Thus
\(\int w\dd\rho_t^\varepsilon\geq B(\rho_0)-D_\varepsilon t\).
Set
\begin{equation}\label{eq:explicit-local-time}
T_0=\min\left\{1,\frac{B(\rho_0)}{2(1+D_\varepsilon)}\right\}>0.
\end{equation}
For $0\leq t\leq T_0$ the denominator is at least
$B(\rho_0)/2>\varepsilon$. Therefore the cutoff is inactive:
\(m_\varepsilon(t) =\frac{\int vw\dd\rho_t^\varepsilon} {\int w\dd\rho_t^\varepsilon}\).
The restriction $\rho=\rho^\varepsilon|_{[0,T_0]}$
is the desired solution of the original nonlinear equation.
Its continuity, initial value, and fourth-moment bound
are inherited from \cref{lem:cutoff-existence}.
Finally \cref{prop:observables}(ii) gives $m\in C^k$
when the higher-order weight condition holds.
\end{proof}

\begin{corollary}[A continuation criterion]
\label[corollary]{cor:existence-continuation}
Assume \eqref{eq:existence-weight}. Let $S<\infty$ and
let $\rho$ be a narrowly continuous probability weak
solution of \eqref{eq:cbo} on $[0,S)$, with initial
value $\rho_0$. If
\(\sup_{0\leq t<S}\int|v|^2\dd\rho_t(v)<\infty\),
then it extends to a solution on $[0,S+\delta]$ for some
$\delta>0$, with $W_2$ continuity and a finite fourth
moment on that interval.
\end{corollary}

\begin{proof}
The moment bound gives $B(\rho_t)\geq\frac12\min_{\overline B_R}w>0$
for some fixed $R$, so $m$ is bounded.
By \cref{lem:fourth-moment} and the weak equation,
$\rho_t\to\rho_S\in\Pp_4$ in $W_2$ as $t\uparrow S$.
Restarting from $\rho_S$, \eqref{eq:explicit-local-time} gives
$\delta=\min\{1,B(\rho_S)/[2(1+D_{B(\rho_S)/4})]\}>0$.
Joining this solution to $\rho$ at $S$ gives the extension to $[0,S+\delta]$.
\end{proof}

\subsection{Global existence}
\label{sec:linear-representation}

\begin{theorem}[Global well-posedness and continuation]\label{thm:global-existence}
Under \eqref{eq:existence-weight}, every local weak solution of
\eqref{eq:cbo} has a unique continuation to arbitrary finite times.
In particular, there is a unique global weak solution with initial
value $\rho_0$. It belongs to $C([0,T];\Pp_2)$ in
the $W_2$ topology and has a bounded fourth moment on every finite
interval $[0,T]$.
\end{theorem}

For a prescribed continuous trajectory $m$, let $W^{ij}$ be $d^2$
independent standard Brownian motions and solve
\begin{equation}\label{eq:global-sde}
\dd X_s=-\lambda(X_s-m(s))\dd s
 +\sigma\sum_{i,j=1}^d(X_{s,j}-m_j(s))e_i\dd W_s^{ij}.
\end{equation}
Its covariance is $\sigma^2|X_s-m(s)|^2I_d$, the same as
for the usual noise $\sigma|X_s-m(s)|\dd B_s$ with $d$-dimensional
Brownian motion $B$. Thus its generator is \eqref{eq:generator}.
We use \eqref{eq:global-sde} because its coefficients are affine
in space and globally Lipschitz, giving a unique strong solution.

By \cref{lem:fourth-moment}, on every compact time interval,
$\int_0^T\int(|v-m(t)|^2+|v-m(t)|^4)\dd\rho_t(v)\dd t<\infty$.
Thus \cite[Theorem~1.3, p.~3232]{LackerShkolnikovZhang2023}, applied
with deterministic coefficients, zero common noise, and $p=2$, represents
every weak solution with prescribed $m$ as the marginals of a martingale
solution with generator $L_{m(t)}$. The globally Lipschitz coefficients of
\eqref{eq:global-sde} give uniqueness for this martingale problem, so these
marginals coincide with its SDE law.

It\^o's formula gives the second-moment bound
\begin{equation}\label{eq:global-moment}
\mathbb E|X_t|^2\leq C_T\left(\mathbb E|X_0|^2+
                              \int_0^t|m(s)|^2\dd s\right),
\qquad t\leq T.
\end{equation}
Indeed its exact derivative is
$(-2\lambda+d\sigma^2)\mathbb E|X_t|^2+
2(\lambda-d\sigma^2)m(t)\cdot\mathbb E X_t+d\sigma^2|m(t)|^2$.
Young's inequality and Gronwall give \eqref{eq:global-moment}, with
$C_T$ depending only on $d,\lambda,\sigma,T$, not on $m$.

\begin{lemma}[A lower bound for the transition law]\label[lemma]{lem:global-lower}
There are constants $h_0,\varepsilon_0,c_0>0$, depending only on
$d,\lambda,\sigma$, with the following property. Let
$0<t-a\leq h_0$, let $m\in C([a,t];\R^d)$, and put
\(Q=\int_a^t|m(s)|^2\dd s>0.\)
For the solution of \eqref{eq:global-sde} with $X_a=x$, assume
$|x|/\sqrt Q\leq\varepsilon_0$. Then its terminal law satisfies
\(\mathbb P(X_t\in E)\geq c_0Q^{-d/2}
                  |E\cap B_{\sigma\sqrt Q/2}(0)|\)
for every Borel set $E$. This is a lower bound as measures;
no density theorem is presumed. The constants are independent of
$a,t,Q,x$ and $m$.
\end{lemma}
\begin{remark}AI tools assisted with the proof of this lemma, particularly in
identifying the result in
\cite[Proposition~3.1(i)]{BallyCaramellino2013}.
For a prescribed trajectory $m(t)$, the SDE has affine coefficients,
so the Malliavin derivative formulas in \cite{Nualart2006} apply
directly, making the required estimates straightforward.
\end{remark}
\begin{proof}
Put $Y_s=X_s/\sqrt Q$ and $n(s)=m(s)/\sqrt Q$.
Then $\int_a^t|n|^2=1$, $\int_a^t|n|\leq\sqrt{t-a}$, and
\begin{equation}\label{eq:global-Y}
\dd Y_s=-\lambda(Y_s-n(s))\dd s
 +\sigma\sum_{i,j}(Y_{s,j}-n_j(s))e_i\dd W_s^{ij},
\qquad Y_a=x/\sqrt Q.
\end{equation}
Separating the integral containing $n$ from the stochastic integral gives
\begin{equation}\label{eq:global-gaussian-remainder}
\begin{aligned}
Y_t&=x/\sqrt Q-\sigma Z+R_t,
& Z_i&=\sum_j\int_a^t n_j(s)\dd W_s^{ij},\\
R_t&=-\lambda\int_a^tY_s\dd s+\lambda\int_a^t n(s)\dd s
 +\sigma\sum_{i,j}\int_a^tY_{s,j}e_i\dd W_s^{ij}.
\end{aligned}
\end{equation}
The independent Brownian rows and $\int_a^t|n|^2=1$ imply
$Z\sim N(0,I_d)$. Thus \eqref{eq:global-gaussian-remainder} has the form
$F=x+G+R$ of \cite[Section~3.1, Proposition~3.1(i)]{BallyCaramellino2013},
with shift $x/\sqrt Q$, Gaussian term $G=-\sigma Z$, and remainder $R_t$.
Here $G$ is an It\^o integral with deterministic integrands and covariance
$\sigma^2I_d$. To apply \cite[Proposition~3.1(i)]{BallyCaramellino2013},
we need to establish
$\|R_t\|_{\mathbb D^{2,p}}\leq C_p\sqrt{t-a}(1+|x|/\sqrt Q)$
for $p\geq2$, where
$\|F\|_{\mathbb D^{2,p}}=\|F\|_p+
\|DF\|_{L^p(\mathcal H)}+\|D^2F\|_{L^p(\mathcal H^{\otimes2})}$.
Here $\mathcal H=L^2([a,t];\R^{d^2})$ and $D_r^{k\ell}$ denotes
Malliavin differentiation at time $r$ with respect to $W^{k\ell}$.
The $L^p$ norms are over the probability space:
$\|F\|_p=(\mathbb E|F|^p)^{1/p}$,
$\|DF\|_{L^p(\mathcal H)}=(\mathbb E\|DF\|_{\mathcal H}^p)^{1/p}$,
and $\|D^2F\|_{L^p(\mathcal H^{\otimes2})}
=(\mathbb E\|D^2F\|_{\mathcal H^{\otimes2}}^p)^{1/p}$.
Throughout, $t-a\leq1$, $p\geq2$, and $C_p$ depends only on $d,\lambda,\sigma,p$.

We first estimate $Y$ and its first two Malliavin derivatives, then
substitute these bounds into the expression for $R_t$.
Since $\lambda\int_a^t|n|\leq\lambda\sqrt{t-a}$ and
$\int_a^t|n|^2=1$, the Burkholder--Davis--Gundy inequality gives
\[
\begin{aligned}
\mathbb E\sup_{a\leq r\leq s}|Y_r|^p
&\leq C_p\left(1+\frac{|x|^p}{Q^{p/2}}\right)
 +C_p\mathbb E\left(\int_a^s|Y_u|\dd u\right)^p
 +C_p\mathbb E\left(\int_a^s|Y_u|^2\dd u\right)^{p/2}\\
&\leq C_p\left(1+\frac{|x|^p}{Q^{p/2}}\right)
 +C_p\int_a^s\mathbb E\sup_{a\leq r\leq u}|Y_r|^p\dd u.
\end{aligned}
\]
The second inequality follows from
$(\int_a^s|Y_u|\dd u)^p\leq(s-a)^{p-1}\int_a^s|Y_u|^p\dd u$,
$(\int_a^s|Y_u|^2\dd u)^{p/2}
\leq(s-a)^{p/2-1}\int_a^s|Y_u|^p\dd u$, and $s-a\leq1$.
Gronwall's inequality yields
$\|\sup_{a\leq s\leq t}|Y_s|\|_p\leq C_p(1+|x|/\sqrt Q)$.

Apply \cref{lem:affine-malliavin} below to \eqref{eq:global-Y}.
The first term in the equation for $D_r^{k\ell}Y_s$ is
$\sigma e_k(Y_{r,\ell}-n_\ell(r))$, with $D_rY_s=0$ for $r>s$.
Its squared $\mathcal H$ norm satisfies
\[
\sigma^2\sum_{k,\ell}\int_a^s
 |e_k(Y_{r,\ell}-n_\ell(r))|^2\dd r
 =d\sigma^2\int_a^s|Y_r-n(r)|^2\dd r
 \leq2d\sigma^2\left((t-a)\sup_{a\leq r\leq t}|Y_r|^2+1\right).
\]
Taking the $p/2$ power and expectation, and then using the bound on $Y$,
we obtain
\[
\begin{aligned}
\mathbb E\left(\int_a^s|Y_r-n(r)|^2\dd r\right)^{p/2}
&\leq C_p\left((t-a)^{p/2}
       \mathbb E\sup_{a\leq r\leq t}|Y_r|^p+1\right)\\
&\leq C_p\left((t-a)^{p/2}
       \left(1+\frac{|x|^p}{Q^{p/2}}\right)+1\right)
 \leq C_p\left(1+\frac{|x|^p}{Q^{p/2}}\right).
\end{aligned}
\]
The stochastic integrand satisfies
$\sigma^2\sum_{i,j}\|(DY_u)_j e_i\|_{\mathcal H}^2
=d\sigma^2\|DY_u\|_{\mathcal H}^2$.
The Hilbert-space Burkholder--Davis--Gundy inequality bounds the
stochastic integral, while H\"older's inequality bounds the drift.
Combining these integral estimates with the preceding bound for the first term gives
\[
\begin{aligned}
\mathbb E\sup_{a\leq v\leq s}\|DY_v\|_{\mathcal H}^p
&\leq C_p\mathbb E\left(\int_a^s|Y_r-n(r)|^2\dd r\right)^{p/2}\\
&\quad+C_p\mathbb E\left(\int_a^s\|DY_u\|_{\mathcal H}\dd u\right)^p
 +C_p\mathbb E\left(\int_a^s\|DY_u\|_{\mathcal H}^2\dd u\right)^{p/2}\\
&\leq C_p\left(1+\frac{|x|^p}{Q^{p/2}}\right)
 +C_p\int_a^s\mathbb E\sup_{a\leq v\leq u}
                       \|DY_v\|_{\mathcal H}^p\dd u.
\end{aligned}
\]
Gronwall's inequality yields
$\|\sup_{a\leq s\leq t}\|DY_s\|_{\mathcal H}\|_p
\leq C_p(1+|x|/\sqrt Q)$.

For $r<q\leq s$, the first term in the second derivative equation
of \cref{lem:affine-malliavin} is $\sigma e_u(D_r^{k\ell}Y_q)_v$.
For $q<r$, interchange $(q,u,v)$ and $(r,k,\ell)$.
The first terms for $r<q$ and $q<r$ have combined squared norm
\[
\begin{aligned}
&\sigma^2\sum_{k,\ell,u,v}\int_a^s\int_a^q
       |e_u(D_r^{k\ell}Y_q)_v|^2\dd r\dd q\\
&\quad+\sigma^2\sum_{k,\ell,u,v}\int_a^s\int_a^r
       |e_k(D_q^{uv}Y_r)_\ell|^2\dd q\dd r
 =2d\sigma^2\int_a^s\|DY_q\|_{\mathcal H}^2\dd q.
\end{aligned}
\]
Thus the bound on $DY$ gives
\[
\mathbb E\left(2d\sigma^2\int_a^s\|DY_q\|_{\mathcal H}^2\dd q\right)^{p/2}
\leq C_p(t-a)^{p/2}
           \mathbb E\sup_{a\leq q\leq t}\|DY_q\|_{\mathcal H}^p
\leq C_p(t-a)^{p/2}\left(1+\frac{|x|^p}{Q^{p/2}}\right).
\]
Applying the same Hilbert-space Burkholder--Davis--Gundy and H\"older
inequalities to the second derivative equation now gives
\[
\mathbb E\sup_{a\leq v\leq s}\|D^2Y_v\|_{\mathcal H^{\otimes2}}^p
\leq C_p(t-a)^{p/2}\left(1+\frac{|x|^p}{Q^{p/2}}\right)
 +C_p\int_a^s\mathbb E\sup_{a\leq v\leq u}
                \|D^2Y_v\|_{\mathcal H^{\otimes2}}^p\dd u.
\]
Gronwall's inequality yields
$\|\sup_{a\leq s\leq t}\|D^2Y_s\|_{\mathcal H^{\otimes2}}\|_p
\leq C_p\sqrt{t-a}(1+|x|/\sqrt Q)$.

We can now estimate the three terms defining
$\|R_t\|_{\mathbb D^{2,p}}$.
For $R_t$ itself, Minkowski's and the Burkholder--Davis--Gundy
inequalities in \eqref{eq:global-gaussian-remainder} give
\[
\|R_t\|_p \leq\lambda(t-a)\left\|\sup_{a\leq s\leq t}|Y_s|\right\|_p +\lambda\sqrt{t-a} +C_p\left\|\left(\int_a^t|Y_s|^2\dd s\right)^{1/2}\right\|_p \leq C_p\sqrt{t-a}\left(1+\frac{|x|}{\sqrt Q}\right).
\]
Since $D_r^{k\ell}Z_i=\delta_{ik}n_\ell(r)$,
\eqref{eq:global-gaussian-remainder} gives
$D_r^{k\ell}R_t=D_r^{k\ell}Y_t+\sigma e_kn_\ell(r)$.
In the first derivative equation of \cref{lem:affine-malliavin},
this replaces $\sigma e_k(Y_{r,\ell}-n_\ell(r))$ by
$\sigma e_kY_{r,\ell}$ and leaves both integrals unchanged.
Applying Minkowski's and the Burkholder--Davis--Gundy inequalities,
then substituting the bounds on $Y$ and $DY$, gives
\[
\begin{aligned}
\|DR_t\|_{L^p(\mathcal H)}
&\leq\sigma\sqrt d\left\|\left(\int_a^t|Y_r|^2\dd r\right)^{1/2}\right\|_p
 +\lambda\int_a^t\|DY_s\|_{L^p(\mathcal H)}\dd s
 +C_p\left\|\left(\int_a^t\|DY_s\|_{\mathcal H}^2\dd s\right)^{1/2}\right\|_p\\
&\leq C_p\sqrt{t-a}\left\|\sup_s|Y_s|\right\|_p
 +C_p\big((t-a)+\sqrt{t-a}\big)
           \left\|\sup_s\|DY_s\|_{\mathcal H}\right\|_p\\
&\leq C_p\sqrt{t-a}\left(1+\frac{|x|}{\sqrt Q}\right).
\end{aligned}
\]
Finally, $D^2Z=0$, so \eqref{eq:global-gaussian-remainder}
gives $D^2R_t=D^2Y_t$.
The bound on $D^2Y$ therefore gives
$\|D^2R_t\|_{L^p(\mathcal H^{\otimes2})}
\leq C_p\sqrt{t-a}(1+|x|/\sqrt Q)$.
Adding the three estimates proves
\begin{equation}\label{eq:global-remainder-malliavin}
\|R_t\|_{\mathbb D^{2,p}}
 \leq C_p\sqrt{t-a}\left(1+\frac{|x|}{\sqrt Q}\right).
\end{equation}

We apply \cite[Proposition~3.1(i)]{BallyCaramellino2013} with this remainder.
It supplies a random variable $0\leq U\leq1$ on the underlying
probability space and a nonnegative function $\varrho_{Y_t,U}$ such that
$\mathbb E[U\mathbf1_{\{Y_t\in E\}}]=\int_E \varrho_{Y_t,U}(y)\dd y$
for every Borel set $E\subset\R^d$, with
$\sup_y|\varrho_{Y_t,U}(y)-\gamma_{\sigma^2I_d}(y-x/\sqrt Q)|
\leq c_d\sigma^{-d}(1+\|R_t/\sigma\|_{\mathbb D^{2,q_d}})^{\ell_d}
\|R_t/\sigma\|_{\mathbb D^{2,q_d}}$.
Here $\gamma_{\sigma^2I_d}$ is the centered Gaussian density with covariance
$\sigma^2I_d$, and $c_d,q_d,\ell_d$ depend only on $d$.
Choose $\varepsilon_0=\min\{1,\sigma/2\}$.
For $|y|\leq\sigma/2$ and $|x|/\sqrt Q\leq\varepsilon_0$,
$|y-x/\sqrt Q|\leq\sigma$, hence
$\gamma_{\sigma^2I_d}(y-x/\sqrt Q)\geq(2\pi\sigma^2)^{-d/2}e^{-1/2}$.
Fix $p\geq\max\{2,q_d\}$ in terms of $d$.
By \eqref{eq:global-remainder-malliavin},
$\|R_t/\sigma\|_{\mathbb D^{2,q_d}}\leq(2C_p/\sigma)\sqrt{t-a}$.
We choose $h_0\leq1$ so that $(2C_p/\sigma)\sqrt{h_0}\leq1$ and
$c_d\sigma^{-d}2^{\ell_d}(2C_p/\sigma)\sqrt{h_0}
\leq\tfrac12(2\pi\sigma^2)^{-d/2}e^{-1/2}$.
Then $h_0$ depends only on $d,\lambda,\sigma$, and the density comparison
gives $\varrho_{Y_t,U}(y)\geq\tfrac12(2\pi\sigma^2)^{-d/2}e^{-1/2}$
on $B_{\sigma/2}$.
Since $0\leq U\leq1$, every Borel set $E$ satisfies
\[
\mathbb P(Y_t\in E) \geq\mathbb E[U\mathbf1_{\{Y_t\in E\cap B_{\sigma/2}\}}] =\int_{E\cap B_{\sigma/2}}\varrho_{Y_t,U}(y)\dd y \geq\tfrac12(2\pi\sigma^2)^{-d/2}e^{-1/2}|E\cap B_{\sigma/2}|.
\]
Applying this to $E/\sqrt Q$ and using
$|(E/\sqrt Q)\cap B_{\sigma/2}|
=Q^{-d/2}|E\cap B_{\sigma\sqrt Q/2}|$ proves the lemma with
$c_0=\tfrac12(2\pi\sigma^2)^{-d/2}e^{-1/2}$.
\end{proof}

\begin{lemma}[Malliavin derivatives of affine SDEs]
\label[lemma]{lem:affine-malliavin}
Let $B,A_1,\ldots,A_N\in C([a,T];\R^{d\times d})$ and
$b,c_1,\ldots,c_N\in C([a,T];\R^d)$ be deterministic.
Let $W^1,\ldots,W^N$ be independent standard Brownian motions, and
let $X$ solve
$\dd X_s=(B(s)X_s+b(s))\dd s+
\sum_{\alpha=1}^N(A_\alpha(s)X_s+c_\alpha(s))\dd W_s^\alpha$
with deterministic $X_a=x\in\R^d$.
Then $X_s\in\mathbb D^{2,p}$ for every $p\geq2$.
Writing $D_r^\alpha$ for Malliavin differentiation with respect to
$W^\alpha$ at time $r$, the first derivative satisfies, for $a\leq r\leq s$,
\[
D_r^\alpha X_s=A_\alpha(r)X_r+c_\alpha(r)
 +\int_r^s B(u)D_r^\alpha X_u\dd u
 +\sum_{\gamma=1}^N\int_r^s A_\gamma(u)D_r^\alpha X_u\dd W_u^\gamma,
\]
with $D_r^\alpha X_s=0$ for $r>s$.
For $a\leq r<q\leq s$, the second derivative satisfies
\[
D_q^\beta D_r^\alpha X_s=A_\beta(q)D_r^\alpha X_q
 +\int_q^s B(u)D_q^\beta D_r^\alpha X_u\dd u
 +\sum_{\gamma=1}^N\int_q^s A_\gamma(u)D_q^\beta D_r^\alpha X_u\dd W_u^\gamma.
\]
For $q<r$, interchange $(q,\beta)$ and $(r,\alpha)$;
the second derivative is zero when $s<\max\{r,q\}$.
\end{lemma}

\begin{proof}
The coefficients have bounded spatial derivatives of positive order
and bounded values at the origin on $[a,T]$.
Thus \cite[Theorem~2.2.2, pp.~122--124]{Nualart2006} gives
$X_s\in\mathbb D^{2,p}$.
The spatial Jacobians are $B(s)$ and $A_\alpha(s)$, so
\cite[Theorem~2.2.1, Eq.~(2.41), p.~119]{Nualart2006} gives the first equation.
For $r<q$, $D_q^\beta(A_\alpha(r)X_r+c_\alpha(r))=0$.
Applying \cite[Lemma~2.2.2, Eq.~(2.52), p.~122]{Nualart2006}
to the first equation gives the second: the stochastic integral
contributes $A_\beta(q)D_r^\alpha X_q$, while the deterministic matrices
have zero Malliavin derivatives.
Symmetry gives the case $q<r$.
\end{proof}

\begin{proof}[Proof of Theorem~\ref{thm:global-existence}]
Let $S$ be the endpoint of a maximal continuation of a given local
solution, and suppose that $S<\infty$.

If $\int_0^S|m(s)|^2\dd s<\infty$, the moment estimate
\eqref{eq:global-moment} bounds the second moment uniformly for $t<S$.
\Cref{cor:existence-continuation} contradicts maximality.
It remains to exclude divergence of this integral.

Choose $a<S$ with $S-a\leq h_0$ and fix the marginal $\rho_a$.
Choose a finite radius $R$ with $\rho_a(B_R)\geq1/2$.
For $a<t<S$ now put
$Q(t)=\int_a^t|m(s)|^2\dd s$.
If the total integral diverges, then $Q(t)\to\infty$.
For every sufficiently late $t$, $Q(t)\geq R^2/\varepsilon_0^2$,
so every $x\in B_R$ satisfies $|x|/\sqrt{Q(t)}\leq\varepsilon_0$.
The linear representation on $[a,t]$ starts with distribution $\rho_a$.
Its marginal at $t$ is the mixture of the transition laws started
at the points $x$ distributed according to $\rho_a$. Since
$t-a<S-a\leq h_0$, Lemma~\ref{lem:global-lower} applies to every
such starting point in $B_R$ and gives
\[
\begin{aligned}
 \rho_t(E)
 &=\int_{\R^d}\mathbb P(X_t\in E\mid X_a=x)\dd\rho_a(x)\\
 &\geq\int_{B_R}c_0Q(t)^{-d/2}
                 |E\cap B_{\sigma\sqrt{Q(t)}/2}|\dd\rho_a(x)\\
 &\geq\frac{c_0}{2}Q(t)^{-d/2}
                 |E\cap B_{\sigma\sqrt{Q(t)}/2}|.
\end{aligned}
\]
Here the conditional law denotes the solution of the prescribed
linear SDE started at $x$ at time $a$; the last inequality uses
$\rho_a(B_R)\geq1/2$. This is an inequality between measures, so it
also holds after integrating any nonnegative measurable function.
For all sufficiently late $t$ we also have
$Q(t)\geq\max\{1,4/\sigma^2\}$, which ensures
$B_1\subset B_{\sigma\sqrt{Q(t)}/2}$. Applying the measure inequality
to $w\mathbf1_{B_1}$ gives
\(\int w\dd\rho_t \geq\int_{B_1}w\dd\rho_t \geq\frac{c_0}{2}Q(t)^{-d/2}\int_{B_1}w(v)\dd v.\)
Consequently,
\begin{equation}\label{eq:global-B-lower}
B(t):=\int w\dd\rho_t\geq
 \frac{c_0}{2}\left(\int_{B_1}w(v)\dd v\right)Q(t)^{-d/2}.
\end{equation}
The prefactor is strictly positive and depends only on
$d,\lambda,\sigma$ and $w$. The time after which the estimate
holds may additionally depend on $\rho_a$ and the trajectory.
To estimate the consensus, split its numerator at radius
$\sqrt{Q(t)}$. On the inner ball, positivity of $w$ gives
\[
 \int_{|v|\leq\sqrt{Q(t)}}|v|w(v)\dd\rho_t(v)
 \leq\sqrt{Q(t)}\int_{|v|\leq\sqrt{Q(t)}}w(v)\dd\rho_t(v)
 \leq\sqrt{Q(t)}B(t).
\]
On the complement, use that $\rho_t$ is a probability measure:
\[
\int_{|v|>\sqrt{Q(t)}}|v|w(v)\dd\rho_t(v) \leq\left(\sup_{|v|>\sqrt{Q(t)}}|v|w(v)\right) \rho_t(\{|v|>\sqrt{Q(t)}\}) \leq\sup_{|v|>\sqrt{Q(t)}}|v|w(v).
\]
Adding the two bounds and dividing by $B(t)>0$ proves the first
two inequalities below:
\begin{equation}\label{eq:global-m-bound}
|m(t)| \leq\frac{\int |v|w(v)\dd\rho_t(v)}{B(t)} \leq\sqrt{Q(t)} +\frac{\sup_{|v|>\sqrt{Q(t)}}|v|w(v)}{B(t)} \leq\sqrt{Q(t)}+C.
\end{equation}
For the last inequality, rapid decay of $w$ ensures that
$C_w:=\sup_{v\in\R^d}(1+|v|)^{d+2}w(v)<\infty$.
For $|v|>\sqrt Q\geq1$ we have
$|v|w(v)\leq C_w|v|^{-(d+1)}\leq C_wQ^{-(d+1)/2}$.
Combining this bound with \eqref{eq:global-B-lower} gives explicitly
\[
 \frac{\sup_{|v|>\sqrt{Q(t)}}|v|w(v)}{B(t)}
 \leq\frac{2C_w}{c_0\int_{B_1}w(v)\dd v}
             Q(t)^{-(d+1)/2}Q(t)^{d/2}
 \leq\frac{2C_w}{c_0\int_{B_1}w(v)\dd v}.
\]
One may therefore take $C=2C_w/(c_0\int_{B_1}w(v)\dd v)$ in
\eqref{eq:global-m-bound}. Thus the constant in \eqref{eq:global-m-bound} depends only on $d,\lambda,\sigma$ and the weight, through $\int_{B_1}w$ and its decay bound.

Since $Q'(t)=|m(t)|^2$, \eqref{eq:global-m-bound} yields
\(Q'(t)\leq2Q(t)+2C^2\)
for all sufficiently late $t<S$. Choose one such $t_0<S$.
Gronwall gives
$Q(t)\leq e^{2(t-t_0)}Q(t_0)+C^2(e^{2(t-t_0)}-1)$,
which stays finite as $t\uparrow S$. This contradicts the
divergence of $Q$. Both alternatives for the forcing integral
are excluded, so there is no finite maximal endpoint.

Finally, on each finite closed interval narrow continuity makes
$B(t)$ continuous and positive, hence bounded below. The consensus
is bounded, the established fourth-moment estimate applies,
and narrow continuity together with uniform integrability of
the second moment gives $W_2$ continuity. Uniqueness follows from
\cref{prop:weak-uniqueness} below.
\end{proof}

\subsection{Uniqueness and Wasserstein stability}
\label{sec:weak-uniqueness}

\begin{proposition}[Uniqueness and $W_2$ stability]
\label[proposition]{prop:weak-uniqueness}
Suppose $w>0$ and $w,v_iw\in\mathcal X^2$.
Let $\rho,\nu$ be weak solutions in the sense of
\cref{def:weak-solution} on $[0,T]$, with initial laws
$\rho_0,\nu_0\in\Pp_4$ and the same parameters and weight.
Set $b_T=\min_{0\leq t\leq T}\{B(\rho_t),B(\nu_t)\}>0$ and
$K_T=\|\nabla(vw)\|_\infty/b_T+
\|vw\|_\infty\|\nabla w\|_\infty/b_T^2$.
Then
\begin{equation}\label{eq:weak-W2-stability}
 W_2^2(\rho_t,\nu_t)\leq e^{C_Tt}W_2^2(\rho_0,\nu_0),
 \qquad 0\leq t\leq T,
\end{equation}
where $C_T=-2\lambda+d\sigma^2+
2|\lambda-d\sigma^2|K_T+d\sigma^2K_T^2$.
The constants depend only on $d,\lambda,\sigma$, the displayed
norms of $w,vw$, and the two curves through $b_T$.
In particular, two weak solutions with the same initial law coincide.
\end{proposition}

\begin{proof}
Boundedness and continuity of $w,vw$ make both consensus trajectories
$m_\rho,m_\nu$ continuous. The functions $B(\rho_t),B(\nu_t)$
are continuous and strictly positive on $[0,T]$, which proves $b_T>0$.
\Cref{lem:fourth-moment} applies to each prescribed trajectory, so
both solutions have uniformly bounded fourth moments.
Taking $\varepsilon=b_T$ in
\eqref{eq:cutoff-consensus-lipschitz} gives
$|m_\rho(t)-m_\nu(t)|\leq K_TW_1(\rho_t,\nu_t)
\leq K_TW_2(\rho_t,\nu_t)$.

Choose an optimal quadratic coupling of $\rho_0,\nu_0$ as the joint
law of $X_0,Y_0$, independent of a common Brownian matrix
$(W^{ij})_{i,j=1}^d$. Solve \eqref{eq:global-sde} for $X,Y$ with
the respective prescribed trajectories $m_\rho,m_\nu$ and this same
Brownian matrix. By \cref{sec:linear-representation}, their marginal
laws are the given PDE solutions $\rho_t,\nu_t$.
Put $D(t)=\mathbb E|X_t-Y_t|^2$.
The consensus estimate and this coupling imply
$|m_\rho-m_\nu|\leq K_T\sqrt{D(t)}$.
The diffusion contribution to the derivative of $D$ is
$d\sigma^2\mathbb E|X_t-Y_t-(m_\rho-m_\nu)|^2$.
Thus It\^o's formula gives
\begin{equation}\label{eq:coupled-second-moment}
D'=(-2\lambda+d\sigma^2)D
 +2(\lambda-d\sigma^2)(m_\rho-m_\nu)\cdot\mathbb E(X_t-Y_t)
 +d\sigma^2|m_\rho-m_\nu|^2.
\end{equation}
Since $|\mathbb E(X_t-Y_t)|\leq\sqrt{D(t)}$,
\eqref{eq:coupled-second-moment} gives $D'\leq C_TD$.
Gronwall, $D(0)=W_2^2(\rho_0,\nu_0)$, and
$W_2^2(\rho_t,\nu_t)\leq D(t)$ prove
\eqref{eq:weak-W2-stability}.
\end{proof}

\section{Sharp obstructions at a stationary consensus point}
\label{sec:counterexample}

This appendix proves the two assertions in
\cref{rem:stationary-sharpness}. The calculations were performed with
AI assistance and verified by the author. The counterexamples are
self-consistent solutions of the nonlinear equation
\eqref{eq:cbo}, with a smooth, quadratically coercive
objective and full-support initial densities with all
polynomial moments. The construction uses the Gaussian law of the logarithmic radius;
the singularity analysis uses Fourier analysis.

Take
\begin{equation}\label{eq:counter-data}
f(v)=|v|^2,\qquad
\rho_0(v)=Z_\beta^{-1}(1+|v|^\beta)e^{-|v|^2},
\qquad \beta>-\frac d2,
\end{equation}
where
\(Z_\beta=\int_{\R^d}(1+|v|^\beta)e^{-|v|^2}\dd v\).
For $\beta<0$, the value of the density at $v=0$ does not affect
the probability measure. The condition
$\beta>-d/2$ gives both local integrability and local square
integrability. The density has full support and all
nonnegative polynomial moments. Moreover,
$w=e^{-\alpha|v|^2}$ and $v_iw$ are Schwartz functions, so
\cref{ass:main} holds for every integer $k\geq1$.

We construct a radial solution directly. For any radial
probability density, the oddness of the integrand gives
\(\int_{\R^d}v e^{-\alpha|v|^2}\rho_t(v)\dd v=0\),
and hence $m(t)=0$. The equation to be solved is therefore
\begin{equation}\label{eq:counter-fixed}
\partial_t\rho
 =\lambda\divv_v(v\rho)
       +\frac{\sigma^2}{2}\Delta_v(|v|^2\rho).
\end{equation}
The associated SDE is
$\dd X_t=-\lambda X_t\dd t+\sigma|X_t|\dd W_t$,
with $W$ a $d$-dimensional Brownian motion.
For $X_0\neq0$, It\^o's formula gives
\(\dd\log|X_t| =\left(-\lambda+\frac{(d-2)\sigma^2}{2}\right)\dd t +\sigma\frac{X_t}{|X_t|}\cdot\dd W_t.\)
The last stochastic integral is a one-dimensional Brownian motion,
since its quadratic variation is $t$. Thus
$Z_t:=\log(|X_t|/|X_0|)$ is Gaussian, with
\(\mathbb E Z_t=\left(-\lambda+\frac{(d-2)\sigma^2}{2}\right)t, \operatorname{Var}(Z_t)=\sigma^2t.\)
Rotational invariance gives the density
$\rho_t(v)=\mathbb E[e^{-dZ_t}\rho_0(e^{-Z_t}v)]$.
Substituting \eqref{eq:counter-data} and completing the square yields,
for $t>0$,
\begin{equation}\label{eq:counter-factor}
\rho_t(v) =\frac{e^{d(\lambda+\sigma^2)t}} {Z_\beta\sqrt{2\pi\sigma^2t}} \int_{\R} \exp\!\left( -\frac{\big[z-(\lambda+\frac{(d+2)\sigma^2}{2})t\big]^2} {2\sigma^2t}\right) \times (1+e^{\beta z}|v|^\beta)e^{-e^{2z}|v|^2}\dd z.
\end{equation}
The SDE has globally Lipschitz coefficients, so its continuous paths
give narrow continuity, and It\^o's formula gives
\eqref{eq:counter-fixed} with initial law $\rho_0$.
The radial law has $m(t)=0$ and all polynomial moments, since
$\mathbb E e^{pZ_t}<\infty$ for every $p\geq0$.
Thus it is a global weak solution of \eqref{eq:cbo}.
The positive integrand in \eqref{eq:counter-factor} gives full support.

For each fixed $t>0$, separate the terms $1$ and
$e^{\beta z}|v|^\beta$ in this integral. Both integral
factors, with $|v|^\beta$ placed outside the second
integral, are smooth functions of $v$ on all of
$\R^d$. Indeed, a derivative of order $j$ of
$e^{-e^{2z}|v|^2}$ is bounded by $C_j e^{jz}$,
where $C_j$ depends only on $d$ and $j$. Every such exponential moment
is finite for a Gaussian distribution.
The second smooth factor has value
\begin{equation}\label{eq:counter-cusp-coefficient}
\frac1{Z_\beta}
 \exp\!\left[
 (d+\beta)\left(\lambda+\frac{\sigma^2}{2}(\beta+2)\right)t
 \right]>0
\end{equation}
at $v=0$. Thus the solution is a smooth function plus
$|v|^\beta$ times a smooth factor which does not vanish
near zero. In particular, no cancellation removes the
initial singular power. The same formula shows smoothness
away from the origin at every positive time.

\begin{proposition}[A bounded positive density need not become $C^1$]
\label[proposition]{prop:counter-bounded}
For $0<\beta<1$ in \eqref{eq:counter-data}, the initial
density is bounded, continuous, strictly positive, has
full support and all polynomial moments, but
$\rho_t$ is not $C^1$ in any neighborhood of zero for
every $t>0$. It is smooth away from zero.
\end{proposition}

\begin{proof}
The two smooth factors in \eqref{eq:counter-factor}
depend on $|v|^2$. Their Taylor expansions at zero and
\eqref{eq:counter-cusp-coefficient} give
\(\rho_t(v) =\frac{e^{d(\lambda+\sigma^2)t}}{Z_\beta} + \frac{\exp\!\left[ (d+\beta)(\lambda+\frac{\sigma^2}{2}(\beta+2))t\right]} {Z_\beta}|v|^\beta+O_t(|v|^2)\).
Along $v=re_1$ the difference quotient at zero has a
strictly positive multiple of $r^{\beta-1}$ and is unbounded
as $r\downarrow0$. The assertions about the initial data
follow directly from \eqref{eq:counter-data}.
\end{proof}

\begin{lemma}[Sobolev threshold for a truncated power]
\label[lemma]{lem:counter-power}
Let $\beta>-d/2$, with
$\beta\notin\{0,2,4,\ldots\}$, and let
$\chi\in C_c^\infty(\R^d)$ equal one near zero.
For every $s\geq0$,
\(\chi(v)|v|^\beta\in H^s(\R^d) \quad\Longleftrightarrow\quad s<\beta+\frac d2\).
\end{lemma}

\begin{proof}
As a tempered distribution, $|v|^\beta$ satisfies
\begin{equation}\label{eq:counter-power-fourier}
\widehat{|v|^\beta}(\xi)
 =2^{\beta+d}\pi^{d/2}
   \frac{\Gamma((\beta+d)/2)}{\Gamma(-\beta/2)}
     |\xi|^{-\beta-d},\qquad \xi\neq0.
\end{equation}
Here is a direct derivation. For $-d<\beta<0$, use
\(|v|^\beta =\frac1{\Gamma(-\beta/2)} \int_0^\infty \tau^{-\beta/2-1}e^{-\tau|v|^2}\dd\tau\)
and
$\widehat{e^{-\tau|v|^2}}(\xi)
 =(\pi/\tau)^{d/2}e^{-|\xi|^2/(4\tau)}$.
The substitution $s=|\xi|^2/(4\tau)$ gives
\eqref{eq:counter-power-fourier}.
For $\beta>0$ which is not an even integer, put
$n=\lfloor\beta/2\rfloor$. Integration by parts
$n+1$ times gives
\(|v|^\beta =\frac1{\Gamma(-\beta/2)} \int_0^\infty\tau^{-\beta/2-1} \left[e^{-\tau|v|^2} -\sum_{j=0}^n\frac{(-\tau|v|^2)^j}{j!}\right]\dd\tau\).
This integral converges at both endpoints: the bracket
is $O(\tau^{n+1}|v|^{2n+2})$ at zero, and the polynomial
terms have degree at most $n<\beta/2$ at infinity.
When paired with a Fourier test function supported away
from zero, each polynomial has zero Fourier contribution.
The same Gaussian integral therefore proves
\eqref{eq:counter-power-fourier}. Its coefficient is
nonzero precisely for the values of $\beta$ under
consideration.

The function $(1-\chi(v))|v|^\beta$ is smooth. Its
derivatives of order $N>\beta+d$ are integrable, since
they are $O(|v|^{\beta-N})$ at infinity. On a coordinate
cone $|\xi_j|\geq|\xi|/\sqrt d$, the distributional
Fourier identity
\((i\xi_j)^N\widehat{(1-\chi)|\cdot|^\beta}(\xi) =\widehat{\partial_{v_j}^N((1-\chi)|\cdot|^\beta)}(\xi)\)
then gives arbitrarily rapid decay at large frequencies
by taking $N$ arbitrarily large. Combining this with
\eqref{eq:counter-power-fourier}, we obtain
\(\widehat{\chi|\cdot|^\beta}(\xi) =2^{\beta+d}\pi^{d/2} \frac{\Gamma((\beta+d)/2)}{\Gamma(-\beta/2)} |\xi|^{-\beta-d}+O_N(|\xi|^{-N})\)
for arbitrarily large $N$. The transform is bounded
at low frequencies because $\chi|\cdot|^\beta\in L^1$.
Since the leading coefficient is nonzero, the
high-frequency part of the squared $H^s$ norm has the
same convergence criterion as
\(\int_1^\infty R^{2s-2\beta-d-1}\dd R\).
This is finite exactly when $s<\beta+d/2$ and diverges
logarithmically at equality.
\end{proof}

\begin{theorem}[No universal positive unweighted Sobolev gain]
\label[theorem]{thm:counter-no-gain}
For every $s_0\geq0$ and every $\varepsilon>0$, there
is an initial probability density
$\rho_0\in H^{s_0}(\R^d)$, with full support and all
polynomial moments, for the objective $f(v)=|v|^2$,
whose global weak CBO solution satisfies $m(t)=0$ and
\(\rho_t\in H^{s_0}(\R^d), \rho_t\notin H^{s_0+\varepsilon}_{\mathrm{loc}} \quad\text{near }0, t>0\).
The failure of local membership means that no
neighborhood of zero has this Sobolev regularity.
\end{theorem}

\begin{proof}
Choose $\beta$ such that
\(s_0<\beta+\frac d2<s_0+\varepsilon,  \beta\notin\{0,2,4,\ldots\}\),
and use \eqref{eq:counter-data}. This implies
$\beta>-d/2$. The smooth part of the initial density
is Schwartz. Its singular part near zero is
$|v|^\beta$ times the smooth, nonvanishing factor
$Z_\beta^{-1}e^{-|v|^2}$; away from zero it and all
its derivatives decay rapidly. Multiplication by
a smooth cutoff and \cref{lem:counter-power} therefore
give $\rho_0\in H^{s_0}(\R^d)$.

For every fixed $t>0$, \eqref{eq:counter-factor} writes
the solution near zero as a smooth function plus
$|v|^\beta$ times a smooth factor nonzero at zero,
by \eqref{eq:counter-cusp-coefficient}.
If $\rho_t$ belonged to $H^{s_0+\varepsilon}$ locally
on any neighborhood of zero, subtracting the smooth
term, dividing by this nonvanishing factor, and taking
a smaller compact cutoff would imply
$\chi(v)|v|^\beta\in H^{s_0+\varepsilon}(\R^d)$.
This contradicts \cref{lem:counter-power}.

For global membership in $H^{s_0}$, a Fourier change of variables gives
$\|e^{-dz}\rho_0(e^{-z}\cdot)\|_{H^{s_0}}
\leq e^{-dz/2}\max\{1,e^{-s_0z}\}\|\rho_0\|_{H^{s_0}}$.
Apply Minkowski's inequality to
$\rho_t(v)=\mathbb E[e^{-dZ_t}\rho_0(e^{-Z_t}v)]$ to obtain
$\|\rho_t\|_{H^{s_0}}
\leq\mathbb E[e^{-dZ_t/2}\max\{1,e^{-s_0Z_t}\}]
\|\rho_0\|_{H^{s_0}}<\infty$,
since $Z_t$ is Gaussian.
\end{proof}

\section{Proofs of Lemmas~\ref{lem:fourth-moment} and~\ref{lem:fourier-equation}}
\label{sec:moment-fourier-proofs}

\begin{proof}[Proof of \cref{lem:fourth-moment}]
Put $V(v)=(1+|v|^2)^2$. Its derivatives are
\(\nabla V(v)=4(1+|v|^2)v\), \(\partial_{v_iv_j}V(v) =4(1+|v|^2)\delta_{ij}+8v_iv_j\), and \(\Delta V(v)=4d+4(d+2)|v|^2\).
Substitution into \eqref{eq:generator} gives
\[
L_{m(t)}V(v)
 =-4\lambda(1+|v|^2)\big(|v|^2-m(t)\cdot v\big)
 +2\sigma^2|v-m(t)|^2
                   \big(d+(d+2)|v|^2\big).
\]
Thus $L_{m(t)}V\leq CV$, uniformly on $[0,T]$.
Here and below the constants depend only on $d,\lambda,\sigma$ and $\sup_{[0,T]}|m|$. In the cutoff estimates below they
may additionally depend on $\chi$;
this does not change any of these uniformities.

We cannot yet insert the unbounded function $V$ into
the weak equation, since the fourth moment of $\rho_t$
is the quantity to be proved finite. For $n\geq1$ use
the bounded truncation
\(V_n=\frac{nV}{n+V}\).
Differentiating this expression gives
\(\nabla V_n =\frac{n^2}{(n+V)^2}\nabla V\), and \(\nabla^2V_n =\frac{n^2}{(n+V)^2}\nabla^2V -\frac{2n^2}{(n+V)^3}\nabla V\otimes\nabla V\).
Consequently,
\[
L_{m(t)}V_n
 =\frac{n^2}{(n+V)^2}L_{m(t)}V
   -\frac{\sigma^2n^2}{(n+V)^3}
                    |v-m(t)|^2|\nabla V|^2
 \leq C\frac{n^2V}{(n+V)^2}
   =\frac{Cn}{n+V}V_n
   \leq CV_n.
\]

To justify testing with $V_n$, first subtract its limit
at infinity. Since $V_n-n=-n^2/(n+V)$ and $n+V\geq V$,
the derivative formulas above imply
\(|V_n(v)-n| \leq\frac{n^2}{(1+|v|^2)^2}\), \(|\nabla V_n(v)| \leq\frac{4n^2|v|}{(1+|v|^2)^3} \leq\frac{4n^2}{(1+|v|^2)^{5/2}}\), and \(|\nabla^2V_n(v)| \leq\frac{Cn^2}{(1+|v|^2)^3}\).
For the last estimate, the two terms in the Hessian
formula are bounded using
$|\nabla^2V|\leq C(1+|v|^2)$ and
$|\nabla V|^2=16(1+|v|^2)^2|v|^2
\leq16(1+|v|^2)^3$.
As $|v-m(t)|\leq C(1+|v|^2)^{1/2}$, it follows that
\begin{equation}\label{eq:moment-generator-tail}
|L_{m(t)}(V_n-n)(v)| \leq\lambda|v-m(t)|\,|\nabla V_n(v)| +\frac{\sigma^2}{2}|v-m(t)|^2|\Delta V_n(v)| \leq\frac{Cn^2}{(1+|v|^2)^2}.
\end{equation}

Take the cutoff $\chi_R(v)=\chi(v/R)$ from
\eqref{eq:cutoff}, with $0\leq\chi\leq1$,
$\chi=1$ on $B_1$, and $\supp\chi\subset B_2$.
Since $L_{m(t)}n=0$, we can separate the two errors as
\begin{align*}
L_{m(t)}\big(\chi_R(V_n-n)\big)-L_{m(t)}V_n
&=\Big[
 L_{m(t)}\big(\chi_R(V_n-n)\big)
       -\chi_RL_{m(t)}(V_n-n)\Big]\\
&\qquad\quad +(\chi_R-1)L_{m(t)}V_n.
\end{align*}
Apply \eqref{eq:cutoff} with $\phi=V_n-n$ to the term
in brackets. As in the estimate immediately following
\eqref{eq:cutoff}, it vanishes outside
$R\leq|v|\leq2R$ and on this annulus is bounded by
\(C\big(|V_n-n|+R|\nabla V_n|\big) \leq Cn^2R^{-4}, R\geq1\).
The remaining tail $(\chi_R-1)L_{m(t)}V_n$ is bounded
by $Cn^2R^{-4}$ by \eqref{eq:moment-generator-tail}.
Together with the bound on $V_n-n$, this proves
\[
\|\chi_R(V_n-n)-(V_n-n)\|_\infty+
\sup_{0\leq t\leq T}
 \big\|L_{m(t)}\big(\chi_R(V_n-n)\big)-L_{m(t)}V_n\big\|_\infty
 \leq Cn^2R^{-4}.
\]
Here $C$ depends in addition on $\chi$.

For fixed $n,R$, the function $\chi_R(V_n-n)$ belongs to
$C_c^\infty(\R^d)\subset\mathcal X^2$.
Apply \cref{lem:observable} to this spatial test and integrate
\eqref{eq:observable} from $0$ to $t$. This gives, for every
$t\in[0,T]$,
\(\int\chi_R(V_n-n)\dd\rho_t -\int\chi_R(V_n-n)\dd\rho_0 =\int_0^t\int L_{m(s)}\big(\chi_R(V_n-n)\big)\dd\rho_s\dd s\).

Now keep $n$ fixed and let $R\to\infty$. Because each
$\rho_s$ has mass one, the preceding uniform error bound gives
\[
\left|\int\big[\chi_R(V_n-n)-(V_n-n)\big]\dd\rho_s\right|
\leq\|\chi_R(V_n-n)-(V_n-n)\|_\infty
\leq Cn^2R^{-4}\longrightarrow0,
\]
uniformly for $s\in[0,T]$. This controls both terms on
the left, at $s=t$ and $s=0$. For the time integral on
the right,
\begin{align*}
&\left|\int_0^t\int
 \left[L_{m(s)}\big(\chi_R(V_n-n)\big)-L_{m(s)}V_n\right]
                         \dd\rho_s\dd s\right|\\
&\quad\leq\int_0^t
 \left\|L_{m(s)}\big(\chi_R(V_n-n)\big)-L_{m(s)}V_n\right\|_\infty
                                                       \dd s
 \leq Ct n^2R^{-4}\longrightarrow0.
\end{align*}
Passing to the limit in the compactly supported identity
therefore yields
\(\int(V_n-n)\dd\rho_t =\int(V_n-n)\dd\rho_0 +\int_0^t\int L_{m(s)}V_n\dd\rho_s\dd s\).
Finally, $\int n\dd\rho_t=\int n\dd\rho_0=n$, so the
two constant terms cancel and we obtain
\(\int V_n\dd\rho_t =\int V_n\dd\rho_0+\int_0^t\int L_{m(s)}V_n\dd\rho_s\dd s\).
Since $V_n$ is bounded and continuous, its integral is
continuous in time. The generator inequality now gives
\(\int V_n\dd\rho_t \leq\int V_n\dd\rho_0 +C\int_0^t\int V_n\dd\rho_s\dd s\).
Gronwall's inequality implies
\(\int V_n\dd\rho_t \leq e^{Ct}\int V_n\dd\rho_0 \leq e^{Ct}\int V\dd\rho_0\).
Finally $V_n(v)\uparrow V(v)$ as $n\to\infty$.
Monotone convergence gives the same bound for
$\int V\dd\rho_t$ at every $t\in[0,T]$. Taking the
supremum in time proves \eqref{eq:fourth-moment}.
\end{proof}

\begin{proof}[Proof of \cref{lem:fourier-equation}]
Throughout this proof $C$ depends only on
$d,\lambda,\sigma,T,\|m\|_{C^1([0,T])}$,
$\int(1+|v|^2)^2\dd\rho_0$ and $\chi$.
All dependence on $\xi$ is displayed explicitly.
We first justify the frequency derivatives and their joint
continuity. For $R\geq1$, \eqref{eq:fourth-moment} gives
\[
\sup_{0\leq t\leq T}
 \int_{|v|>R}(1+|v|^2)\dd\rho_t(v)
\leq\frac1{1+R^2}
 \sup_{0\leq t\leq T}\int(1+|v|^2)^2\dd\rho_t(v)
\leq\frac C{1+R^2}\longrightarrow0.
\]
Since $m$ is bounded, the same estimate controls the tails
of integrands bounded by $1+|v-m(t)|^2$. Differentiation
under the integral in \eqref{eq:moving-fourier} is therefore
justified through order two and gives
\(\partial_{\xi_j}F(t,\xi) =-i\int(v_j-m_j(t))e^{-i\xi\cdot(v-m(t))}\dd\rho_t(v)\), and \(\partial_{\xi_i\xi_j}F(t,\xi) =-\int(v_i-m_i(t))(v_j-m_j(t)) e^{-i\xi\cdot(v-m(t))}\dd\rho_t(v)\).
To check joint continuity, first multiply each of these
integrands, and the integrand defining $F$, by the same
spatial cutoff $\chi_R$ as above. For fixed $R$, these
are bounded, compactly supported functions of $v$, varying
continuously in $(t,\xi)$ in the uniform norm on compact
parameter sets. Narrow continuity then makes their integrals
jointly continuous. The tail estimate shows that these
cutoff integrals converge uniformly on compact parameter
sets to $F$ and these derivatives. Thus $F$ and all
its frequency derivatives through order two are jointly
continuous. In particular,
\(\int(v-m(t))e^{-i\xi\cdot(v-m(t))}\dd\rho_t(v) =i\nabla_\xi F(t,\xi)\), and \(\int|v-m(t)|^2e^{-i\xi\cdot(v-m(t))}\dd\rho_t(v) =-\Delta_\xi F(t,\xi)\).

Next fix $\xi,R$. The spatial function
$v\mapsto\chi_R(v)e^{-i\xi\cdot v}$ belongs to
$C_c^\infty(\R^d)\subset\mathcal X^2$.
Applying \cref{lem:observable} to its real and imaginary
parts gives
\(\frac{\dd}{\dd s}\int\chi_R(v)e^{-i\xi\cdot v}\dd\rho_s(v) =\int L_{m(s)}\big[\chi_R(v)e^{-i\xi\cdot v}\big] \dd\rho_s(v)\).
Since $m\in C^1$, multiply this observable by the time
phase $e^{i\xi\cdot m(s)}$ and use the ordinary product
rule. The phase is constant in $v$, so it factors out
of $L_{m(s)}$. Integrating from $0$ to $t$ then yields
\begin{align*}
&\int\chi_R(v)e^{-i\xi\cdot(v-m(t))}\dd\rho_t(v)
 -\int\chi_R(v)e^{-i\xi\cdot(v-m(0))}\dd\rho_0(v)\\
&\quad=\int_0^t\int
  (\partial_s+L_{m(s)})
       \big[\chi_R(v)e^{-i\xi\cdot(v-m(s))}\big]
                                      \dd\rho_s(v)\dd s.
\end{align*}
The three derivatives of the exponential are
\(\partial_s e^{-i\xi\cdot(v-m(s))} =i\,m'(s)\cdot\xi\,e^{-i\xi\cdot(v-m(s))}\), \(\nabla_v e^{-i\xi\cdot(v-m(s))} =-i\xi\,e^{-i\xi\cdot(v-m(s))}\), and \(\Delta_v e^{-i\xi\cdot(v-m(s))} =-|\xi|^2e^{-i\xi\cdot(v-m(s))}\).
Consequently,
\[
(\partial_s+L_{m(s)})e^{-i\xi\cdot(v-m(s))}
 =\left[
 i\,m'(s)\cdot\xi+i\lambda\xi\cdot(v-m(s))
 -\frac{\sigma^2}{2}|\xi|^2|v-m(s)|^2
 \right]e^{-i\xi\cdot(v-m(s))}.
\]

We now remove the spatial cutoff. Applying
\eqref{eq:cutoff} to the exponential gives, uniformly in $s$,
\begin{align*}
&\int\left|
 L_{m(s)}\big[\chi_R e^{-i\xi\cdot(v-m(s))}\big]
 -\chi_RL_{m(s)}e^{-i\xi\cdot(v-m(s))}
 \right|\dd\rho_s(v)\\
&\quad\leq C(1+R|\xi|)\rho_s(\{|v|\geq R\})
 \leq C\big(R^{-4}+|\xi|R^{-3}\big).
\end{align*}
Here the commutator is supported on $R\leq|v|\leq2R$,
the exponential has absolute value one and spatial gradient
of absolute value $|\xi|$, and
$\rho_s(\{|v|\geq R\})\leq R^{-4}\int|v|^4\dd\rho_s$
by \eqref{eq:fourth-moment}. The endpoint cutoff errors
are at most $CR^{-4}$. The remaining error in the time
integrand is
$(\chi_R-1)(\partial_s+L_{m(s)})e^{-i\xi\cdot(v-m(s))}$.
Its absolute integral is bounded by
\(C(|\xi|+|\xi|^2) \int_{|v|>R}(1+|v|^2)\dd\rho_s(v)\),
which tends to zero by the tail estimate above. These
estimates are uniform in $s\in[0,T]$ and locally uniform
in $\xi$, so they also justify the limit under the time
integral.

Letting $R\to\infty$ and substituting the frequency-derivative
identities gives
\[
F(t,\xi)-F(0,\xi)
 =\int_0^t\Big[
 i\,m'(s)\cdot\xi\,F(s,\xi)
 -\lambda\xi\cdot\nabla_\xi F(s,\xi)
 +\frac{\sigma^2}{2}|\xi|^2\Delta_\xi F(s,\xi)
 \Big]\dd s.
\]
In particular, the drift contributes
$i\lambda\xi\cdot(i\nabla_\xi F)=-\lambda\xi\cdot\nabla_\xi F$,
and the diffusion contributes
$-\frac{\sigma^2}{2}|\xi|^2(-\Delta_\xi F)
=\frac{\sigma^2}{2}|\xi|^2\Delta_\xi F$.
The integrand is jointly continuous, since $m'$, $F$,
$\nabla_\xi F$, and $\Delta_\xi F$ are continuous.
Differentiating this time integral therefore proves
\eqref{eq:fourier-equation}.
Thus $F$ has one continuous time derivative and two
continuous frequency derivatives, all jointly continuous
in $(t,\xi)$, with one-sided time derivatives at the
endpoints. The argument uses the given probability
measures and their moment bound throughout; no density
of $\rho_t$ is presupposed.
\end{proof}

\end{document}